\documentclass[12pt,reqno]{amsart}
\usepackage{lineno} 
\usepackage{mathrsfs}
\usepackage{todonotes}

\usepackage{amsmath,amssymb,latexsym,soul,cite,mathrsfs}
\usepackage{mathtools}
\numberwithin{equation}{section}
\usepackage{xcolor}
\definecolor{meuverde}{RGB}{0,100,0}

\definecolor{meuverd1}{HTML}{006400}

\definecolor{meuverde2}{cmyk}{1,0,1,0.6}

\usepackage{color,enumitem,graphicx}
\usepackage{tikz}
\usepackage{placeins}
\definecolor{RefPurple}{RGB}{88,0,120}

\definecolor{Orange}{RGB}{255,127,0}
\definecolor{Purple}{RGB}{255,127,0}
\usepackage[colorlinks=true,urlcolor=blue,
citecolor=blue,linkcolor=blue,linktocpage,pdfpagelabels,
bookmarksnumbered,bookmarksopen]{hyperref}
\usepackage[english]{babel}
\usepackage{enumitem}

\usepackage[left=1.6cm,right=1.6cm,top=1.7cm,bottom=1.7cm]{geometry}

\makeatletter
\newcommand{\leqnomode}{\tagsleft@true}
\newcommand{\reqnomode}{\tagsleft@false}
\makeatother
\newcommand{\ignore}[1]{}

\date{}
\def\nd{\noindent}
\def\thend{\rule{3mm}{3mm}}

\newtheorem{theorem}{Theorem}[section]
\newtheorem{definition}{Definition}[section]
\newtheorem{cor}[theorem]{Corollary}
\newtheorem{proposition}[theorem]{Proposition}
\newtheorem{lemma}[theorem]{Lemma}
\newtheorem{remark}{Remark}[section]

\renewcommand{\theequation}{\arabic{section}.\arabic{equation}}

\newcommand{\Int}{\displaystyle\int_{\Omega}}

\usepackage[normalem]{ulem}

\thanks{Corresponding author: Claudiney Goulart. E-mail: claudiney@ufj.edu.br}

\begin{document} 
	
	\title[Multiplicity of solutions for a class of logistic problems]{Multiplicity of solutions for a class of logistic problems 
    }
	
	\vspace{1cm}

    \author{M.L.M. Carvalho}
	\address{ M.L.M. Carvalho \newline Universidade Federal de Goi\'as, IME, Goi\^ania-GO, Brazil}
	\email{\tt marcos\_leandro\_carvalho@ufg.br}
	
	\author{C. Goulart}
	\address{C. Goulart \newline Universidade Federal de Jata\'i, Jata\'i-GO, Brazil }
	\email{\tt claudiney@ufj.edu.br}

    \author{M.C.M. Rezende}
	\address{M.C.M. Rezende \newline  Universidade de Brasília, Brasília-DF, Brazil }
	\email{\tt  manuelarezende@unb.br }

\author{R. Ruviaro}
	\address{R. Ruviaro\newline Universidade de Brasília, Brasília-DF, Brazil }
	\email{\tt ruviaro@unb.br}

    \thanks{The first author was partially supported by CNPq with grant 300411/2025-1. 
The third author was partially supported by FAPDF with grant  00193.00001821/2022-21. 
The fourth author was partially supported by CNPq with grant 302745/2025-4.}

	\subjclass[2020]{35J05,35J15,35J20,35J60}
	\keywords{onlinear elliptic equations; logistic problems; Nehari manifold; nonlinear Rayleigh quotient; mountain-pass theorem}

	\begin{abstract}


The present paper investigates the existence and multiplicity of positive solutions
{to} a class of {degenerate} logistic problems by combining the Nehari
manifold {method} with {nonlinear generalized} Rayleigh quotients.
The quotient associated with the Nehari constraint determines the threshold
$\lambda_*$ at which the fibering structure and the Nehari manifold geometry {emerge}, {while} a second
quotient identifies the zero-energy threshold $\lambda^*$ {associated with} the minimizing
solution. Together with {the threshold} $\bar\lambda$
{for the existence of nontrivial solutions}, {these results yield} the strict ordering
\(
0<\lambda_*<\bar\lambda<\lambda^*<\lambda_1(\Omega)<\lambda_1(\Omega_0),
\)
showing that the Nehari geometry becomes available before nontrivial solutions
{exist} and that the zero-energy transition occurs only afterwards. In the
variational range, we obtain a positive minimizing solution on
$\mathcal N_\lambda^+$ and a second positive solution {via} the Mountain Pass
Theorem. {Complementary sub--supersolution and truncation arguments
extend the multiplicity result to the remaining parameter range.}
We also establish {the sign change} of the minimizing energy, the
blow-up and normalized asymptotic profile of the minimizing solutions as
$\lambda\uparrow\lambda_1(\Omega_0)$, and the {convergence of the mountain-pass level}
$c_\lambda^m$ {to zero} {as}
$\lambda\uparrow\lambda_1(\Omega)$.
\end{abstract}
	\maketitle
	\section{Introduction}

    The purpose of this paper is to study the existence, multiplicity and
qualitative behavior of solutions to the following degenerate logistic problem
	\begin{equation}\label{eq1}
		\left\{
            \begin{array}{rcll}\tag{$P_{\lambda}$}
			-\Delta u  &=&  \lambda u-b(x)(|u|^{p-2}u- |u|^{q-2}u)&\mbox{in}   \  \Omega,  \\[2ex]
			u&=&0  &\mbox{on} \ \ \ \ \partial \Omega,
		\end{array}
		\right.
	\end{equation}
    where $\Omega \subset \mathbb{R}^N, N \geq 1,$ is a smooth bounded domain,
$\lambda$ is a real positive parameter,
$2 < q < p < 2^*={2N}/{(N-2)}$, if $N \geq 3$,
and $2^*=\infty$, if $N=1,2$,
and {$b \in L^{\infty}(\overline{\Omega})$ is nonnegative}. 
Moreover, $b(x)=0$ in a connected subset $\Omega_0\subset\subset\Omega$
with positive Lebesgue measure and smooth boundary.

{
The main contribution of this work is a detailed analysis of the
parameter-dependent variational structure of Problem~\eqref{eq1}. {In particular, we identify and compare distinct thresholds arising from the Nehari geometry, establish the existence of nontrivial solutions, and characterize the parameter regimes in which the energy of the minimizing solution changes sign.}  Nonlinear generalized Rayleigh
quotients play a central role in this analysis and, { together with the Nehari manifold method, the Mountain Pass Theorem and comparison arguments, provide the main tools for establishing existence and multiplicity results}.
}

One of the motivations for the present work comes from the seminal contribution by López-Gómez \cite{lopez, LG},  in which the author analyzes how the study of the interplay within the same habitat between the most classical laws of Population Dynamics originates from the modern theory of nonlinear parabolic
differential equations. More specifically, in this study, {the author considers} the following parabolic problem
\begin{equation} \label{lopez}
\left\{
\begin{aligned}
    u_t - \Delta u &= \lambda u - a(x) f(x,u) u && \quad \text{in } \Omega \times (0,\infty), \\[1ex]
    u(x,t) &= 0 && \quad \text{on } \partial\Omega \times (0,\infty), \\[1ex]
    u(x,0) &= u_0(x) && \quad \text{in } \Omega,
\end{aligned}
\right.
\end{equation}
where $\Omega \subset \mathbb{R}^N$ is a bounded domain, {and} $a(\cdot)$ is a nonnegative weight function. Additional conditions are also imposed on the nonlinearity $f(\cdot, \cdot)$ and on the boundary of the domain $\Omega$. 

From the perspective of population dynamics, this class of equations arises as a result of the interaction between two classical models of population growth: the exponential Malthusian model and the logistic Verhulst model. In this context, the function $u(\cdot, \cdot)$ describes the evolution of the density of a species distributed over the inhabited region $\Omega$. The parameter $\lambda$ represents the population growth rate, while $a(\cdot)$ incorporates crowding effects within the subregion $\Omega \setminus \Omega_0$. The area $\Omega_0$ is referred to as the favorable region, since it allows for exponential population growth as proposed by {Malthus' law}.





It is important to refer to the contribution of Fernandes and Maia \cite{fernandesmaia} to the study of this class of problems given by \eqref{lopez}. In this sense {the authors consider the existence and multiplicity of solutions} for the problem
\begin{equation*}
\left\{
\begin{array}{rcll}
\partial_t u &=& \Delta u + \lambda u + b(x)|u|^{\nu - 1}u, & (x,t) \in \Omega \times (0,+\infty), \\[1ex]
u &=& 0, & \text{on } \partial\Omega \times (0,+\infty), \\[1ex]
u(x,0) &=& u_0(x), & x \in \Omega,
\end{array}
\right.
\end{equation*}
where $\Omega \subset \mathbb{R}^N$ ($N \geq 2$) is a bounded smooth domain, $\lambda > 0$, $1 < \nu < 2^* - 1$ and the weight function $b \in C(\overline{\Omega})$ satisfies $b(x) \leq 0$, vanishing in a proper open subset $\Omega_0 \subset\subset \Omega$ with positive measure. Using variational methods and a comparison principle, the authors analyze how the behavior of solutions depends on the initial data in the phase space. They also {explore} the stationary case and {take} a variational approach. In the stationary case, they {ensure} the existence of {a} positive solution and {a} sign-changing solution.

We also refer to the recent work of Cardoso, Furtado and Maia  \cite{furtadomaia}, who studied the existence of positive and sign-changing steady-state solutions to a degenerate logistic problem with a general superlinear nonlinearity:
\begin{equation*}
    \left\{
    \begin{array}{rcll}
    -\Delta u &=& \lambda u - b(x)f(u) & \text{in } \Omega, \\[1ex]
    u &=& 0 & \text{on } \partial\Omega,
    \end{array}
    \right.
\end{equation*}
where $\Omega \subset \mathbb{R}^N$ ($N \geq 2$) is a smooth bounded domain, $\lambda > 0$, and the weight function  $b \in L^\infty(\overline{\Omega})$ satisfies $b(x) \geq 0$ and $b(x) = 0$ on a connected subset $\Omega_0 \subset\subset \Omega$ with positive Lebesgue measure and smooth boundary.  In \cite{furtadomaia} variational techniques, specifically those involving the restriction of the functional to the Nehari manifold, were employed to study how changes in the parameter 
$\lambda$ or the function 
$b$ affect the existence and nonexistence of solutions to the elliptic problem under consideration.

Motivated by the works cited above, we observe that we can combine the Nehari manifold technique with the generalized Rayleigh quotient method to study the existence of multiple solutions for the logistic Problem \eqref{eq1}. Our approach to deal with that problem  is inspired by recent advances in the literature, where nonlinear Rayleigh quotients have been extensively explored. We refer to \cite{yavdat0,yavdat1,yavdat2} as pioneering contributions in this direction. For the interested reader, we also point to  more recent 
contributions \cite{choquard, carvalhojde, carvalhotmna, ed1, goul1, ksilva,KayeAb},  
as well as the references therein. A common interest in these studies is to identify an extremal parameter that ensures the effectiveness of variational techniques such as the Nehari method, by showing that the associated fibering map possesses at least two critical points for every $\lambda$ restricted to a given range, {which} is a strategy that we also adopt in the present paper.

More precisely, the quotient $R_n${, defined in \eqref{Rn}, determines} the threshold $\lambda_*$ at which the Nehari/fibering geometry appears, {whereas} $R_e${, defined in \eqref{Re}, determines} the threshold $\lambda^*$ {at which the energy of the minimizing solution changes sign}. On the other hand, $\bar\lambda$ is the actual threshold for the existence of nontrivial solutions. We {establish} the strict ordering
\(
0<\lambda_*<\bar\lambda<\lambda^*<\lambda_1(\Omega)<\lambda_1(\Omega_0).
\)
Hence, the Nehari geometry becomes available before the first nontrivial solution exists, whereas the zero-energy transition occurs only after the existence threshold has been crossed. {
This separation is a distinctive feature of our approach and provides insights into the variational structure of the problem that cannot be obtained from an existence analysis based solely on sub- and supersolution methods.}

The methods employed {in this work play} complementary roles. The generalized Rayleigh quotients {characterize} the fibering and Nehari geometry, provide the variational threshold information and lead to a minimizing solution $u_\lambda\in\mathcal N_\lambda^+$ in the range where the minimization argument applies. The Mountain Pass Theorem yields a second critical point in the corresponding variational range. {For the remaining parameter interval, we employ sub-supersolution methods combined with truncation arguments.} Thus, these complementary methods complete, rather than replace, the structural information furnished by the nonlinear Rayleigh quotients. The analysis also {establishes} the energy transition at $\lambda^*$, the blow-up of the minimizing solutions as $\lambda\uparrow\lambda_1(\Omega_0)$ and the {convergence} {of the mountain-pass level } $c_\lambda^m$ {to zero} as $\lambda\uparrow\lambda_1(\Omega)$.

From a standard variational point of view, finding weak solutions of the elliptic Problem~\eqref{eq1}
is equivalent to finding critical points of the functional
$E_\lambda:H_0^1(\Omega)\to\mathbb{R}$ given by
\begin{eqnarray}\label{functional}
E_{\lambda}(u)
&=&
\dfrac{1}{2}\Int |\nabla u|^2\,dx
-\dfrac{\lambda}{2}\Int |u|^2\,dx
+\Int b(x)\left(\dfrac{1}{p}|u|^p-\dfrac{1}{q}|u|^q\right)dx,
\qquad u\in H_0^1(\Omega).
\end{eqnarray}

The functional $E_\lambda$ is well defined, and its Gâteaux derivative is given by
\begin{equation}\label{derivada}
E'_{\lambda}(u)\phi
=
\Int \nabla u\cdot\nabla\phi\,dx
-\lambda\Int u\phi\,dx
+\Int b(x)\left(|u|^{p-2}u-|u|^{q-2}u\right)\phi\,dx,
\qquad u,\phi\in H_0^1(\Omega).
\end{equation}

{
Since $p,q>2$, {it follows from the expression for the first derivative in }\eqref{derivada},
for {all} $v,w\in H_0^1(\Omega)$ we have
\[
E_\lambda''(u)(v,w)
=
\int_\Omega \nabla v\cdot\nabla w\,dx
-\lambda\int_\Omega vw\,dx
+(p-1)\int_\Omega b(x)|u|^{p-2}vw\,dx
-(q-1)\int_\Omega b(x)|u|^{q-2}vw\,dx,
\]
and $u\mapsto E_\lambda''(u)$ is continuous. Hence
$E_\lambda\in C^2(H_0^1(\Omega),\mathbb R)$.
}

{ Denoting by $\hat{\Omega}_0:=\{x \in \Omega; ~ b(x)=0\}$, we suppose the following hypothesis:
   \begin{itemize}
   \item [$(b_0)$] $ b \in L^{\infty}(\overline{\Omega}), ~ b \geq 0$ and $\Omega_0=int(\hat{\Omega}_0)\subset \subset \Omega$ is a connected set with smooth boundary, $\Omega_0 \neq \emptyset$ and $|\hat{\Omega}_0\backslash \Omega_0|=0.$ Moreover, we suppose
   {$||b||_{\infty}<\lambda_1(\Omega)/ {C_{p,q}}$, with $C_{p,q}>0$ given by \eqref{const}.}
   \end{itemize}}
   
 
Throughout this work, we denote by $\lambda_1(\Omega)$ { and $\lambda_1(\Omega_0)$} the first eigenvalues of the operator $-\Delta$ in $H^1_0(\Omega)$ and {$H^1_0(\Omega_0)$, respectively. We also denote}  by $\varphi_1>0$ { and $\varphi_1^0>0$} their corresponding eigenfunctions. 

{As an example for a function satisfying $(b_0)$ the reader may consider
\begin{equation*}
b(x)= 
\begin{cases}
\varepsilon \|x\|-\varepsilon, &  1<\|x\|\leq 2, \\
-\varepsilon\|x\|+3\varepsilon, & 2\leq \|x\|<3,\\
0, & \text{otherwise},
\end{cases}
\end{equation*}
for every $x \in B_3(0)=\{x \in \mathbb{R}^N; \|x\|<3\}$ and for $\varepsilon>0$ sufficiently small such that {$||b||_{\infty}<\lambda_1(\Omega)/C_{p,q}$}.}

{The variational analysis below is organized around two generalized Rayleigh quotients. The first one is associated with the Nehari constraint and determines the fibering geometry, whereas the second one is associated with the zero level of the energy and determines its sign along the minimizing solution. These quotients will be introduced after the Nehari manifold and the energy-zero set are defined.}

Using the same ideas introduced in \cite{nehari1,nehari2}, we shall consider the Nehari manifold associated with Problem \eqref{eq1}:
	\begin{equation}
		\begin{array}{rcl}
			\mathcal{N}_{\lambda}
			=\left\{ u\in  H_0^1(\Omega)\setminus \{0\}:||u||^2-\lambda||u||^2_2=-\Int b(x)(|u|^p{-}|u|^q)dx
		\right\}\end{array}.
		\label{nehari}
	\end{equation}
Now, in a similar way to that done by Tarantello \cite{tarantello}, we shall split the Nehari manifold $\mathcal{N}_\lambda$ into three disjoint subsets as follows:
	\begin{eqnarray*}
		\mathcal{N}_\lambda^+&=&\{u\in\mathcal{N}_\lambda: E''_\lambda(u)u^2>0 \}\label{n+};\\ 
		\mathcal{N}_\lambda^-&=&\{u\in\mathcal{N}_\lambda: E''_\lambda(u)u^2<0 \}\label{n-};\\
		\mathcal{N}_\lambda^0&=&\{u\in\mathcal{N}_\lambda: E''_\lambda(u)u^2=0\}\label{n0}.
	\end{eqnarray*}
	  It is important to emphasize that $u \mapsto E_\lambda''(u)u^2$ is well defined for each $u \in H_0^1(\Omega)$.  { More precisely,      
$$E''_{\lambda}(u)u^2 = \|u\|^2 - {\lambda} \|u\|_2^2  + (p-1)\Int{b(x)}|u|^p dx-(q-1)\Int{b(x)}|u|^q dx, \  \  u \in H_0^1(\Omega).$$}
In order to apply the nonlinear generalized Rayleigh quotient, we also need to introduce some definitions. Firstly, we shall consider the  following set, which contains all nontrivial roots of functional $E_\lambda$:
	\begin{eqnarray}\label{conjE}
		\mathcal{E}_{\lambda}
		=\left\{ u \in H_0^1(\Omega) \setminus \{0\}:\dfrac{1}{2}||u||^2-\dfrac{\lambda}{2}||u||^2_2=-\Int{b(x)}\left(\frac{1}{p}|u|^p-\frac{1}{q}|u|^q\right)dx\right\}. 
	\end{eqnarray}
	In view of {the} definitions \eqref{derivada} and \eqref{nehari}, we may note  that 
	\begin{equation}\label{a0i}
		u \in \mathcal{N}_\lambda, \mbox{i.e.}, E'_{\lambda}(u)u=0 \, \, \, \mbox{{if and only if}} \,\, \, \lambda = \dfrac{||u||^2 + \Int b(x)(|u|^p-|u|^q)dx}{\|u\|^2_{2}}.
	\end{equation}
	Similarly, from \eqref{functional} and \eqref{conjE}, we observe that
	\begin{equation}\label{a11i}
		u \in \mathcal{E}_{\lambda}, \mbox{i.e.}, E_{\lambda}(u)=0  \, \, \, \mbox{{if and only if}} \,\, \, \lambda= \dfrac{\frac{1}{2}||u||^2 + \Int b(x)\left(\dfrac{1}{p}|u|^p-\dfrac{1}{q}|u|^q\right)dx}{\frac{1}{2}||u||^2_{2}}.
	\end{equation} 
The main feature of this work is to ensure that the minimization problems
{\begin{equation}\label{mp}
c_\lambda^+ := \inf \{ E_{\lambda}(u) : u \in \mathcal{N}^{+}_{\lambda}\}.
\end{equation}}
and
\begin{eqnarray}
c_\lambda^m := \inf_{\gamma \in \Gamma_\lambda} \max_{r \in [0,1]} E_\lambda(\gamma(r)),
\end{eqnarray}
 where {$\Gamma_\lambda$ denotes the mountain-pass class defined in Section~\ref{SB3}}, are attained by some specific functions.


It is important to emphasize that the sets introduced in \eqref{nehari} and \eqref{conjE} enable the definition of nonlinear generalized Rayleigh quotients, which have been extensively studied in recent years (see \cite{yavdat0,yavdat1,yavdat2}).
A key ingredient in this framework is the set
{
\begin{equation*}\label{D}
    D := \big\{u \in H_0^1(\Omega) \setminus \{0\}: \int_\Omega b(x)|u|^q\,dx>0\big\}.
\end{equation*}
By {assumption $(b_0)$, this set can equivalently be characterized as the collection}  of nonzero functions that do not vanish a.e. on $\Omega\setminus\Omega_0$. }
This set plays a crucial role in ensuring the well-definedness of the Rayleigh quotients.
	Under these conditions, we define the functionals $R_n,R_e:H_0^1(\Omega)\setminus\{0\} \to \mathbb{R}$ 
	associated with the parameter $\lambda>0$ in the following form 
\begin{equation}\label{Rn}
R_n(u) := 
\begin{cases}
\displaystyle \frac{\displaystyle||u||^2 + \int_\Omega b(x)(|u|^p - |u|^q)\,dx}{||u||_2^2}, &  u \in D, \\[2ex]
\displaystyle \frac{||u||^2}{||u||_2^2}, & u \in D^c\setminus\{0\}.
\end{cases}
\end{equation}
and
\begin{equation}\label{Re}
R_e(u) := 
\begin{cases}
\displaystyle \frac{\displaystyle\frac{1}{2}||u||^2 + \int_\Omega b(x)\left(\frac{1}{p}|u|^p - \frac{1}{q}|u|^q\right)\,dx}{\frac{1}{2}||u||_2^2}, &  u \in D, \\[2ex]
\displaystyle \frac{||u||^2}{||u||_2^2}, &  u \in D^c\setminus\{0\}.
\end{cases}
\end{equation}
	Note that if we compare \eqref{a0i} with \eqref{Rn} and \eqref{a11i} with \eqref{Re}, we conclude that 
$$
u \in \mathcal{N}_{\lambda}\cap D \text{ {if and only if} } R_n(u) = \lambda
\text{ and }
u \in \mathcal{E}_{\lambda}\cap D \text{ {if and only if} } R_e(u) = \lambda.
$$

At this point it is important to observe that $R_e$ is used in order to decide the sign for solution obtained as a minimizer for the functional $E_\lambda$ restricted to $\mathcal{N}_\lambda^+$ (see Theorem \ref{t2}).
	
By using \eqref{Rn} and \eqref{Re} we obtain  
\begin{equation}\label{Rn,t}
R_n(tu)=\dfrac{||u||^2}{||u||^2_2}+\dfrac{1}{||u||^2_2}\Int b(x)(t^{p-2}|u|^p-t^{q-2}|u|^q)dx, \ \ u \in D,\ \ t\geq  0,
\end{equation}
and  
\begin{equation}\label{Re,t}
R_e(tu)=\dfrac{||u||^2}{||u||^2_2}+\dfrac{2}{||u||^2_2}\Int b(x)\left(\dfrac{t^{p-2}}{p}|u|^p-\dfrac{t^{q-2}}{q}|u|^q\right)dx, \ \ u \in  D,\ \ t\geq 0.
\end{equation}
At this stage we define
\begin{equation}\label{lambdane}
\lambda_*:=\inf_{u\in D}\inf_{t>0}R_n(tu), \;\; \lambda^*:=\inf_{u\in D}\inf_{t>0}R_e(tu)\;\;\text{and}\;\;  {\bar{\lambda}:=\inf\{\lambda>0 :~\eqref{eq1}~\mbox{has a nontrivial weak solution}\}}.
\end{equation}

In the following we will state the main results of this paper. Theorem \ref{t2} provides the existence of a solution to Problem \eqref{eq1}, for each {$\overline{\lambda}<\lambda<\lambda_1(\Omega_0)$}. {The multiplicity of solutions, for each $\lambda \in {\color{blue}(\bar\lambda,\lambda_1(\Omega))}$, is established by Theorem \ref{t1}.} Theorem \ref{t3} ensures the existence of at least one solution when $\lambda=\overline{\lambda}$, while Theorem \ref{t4} guarantees the {nonexistence of nontrivial solutions} to Problem \eqref{eq1} when $\lambda<\overline{\lambda}$ {and the nonexistence of positive solutions when} $\lambda\geq \lambda_1(\Omega_0)$.

\begin{theorem}\label{t2}
Suppose that $2<q<p<2^*$ and  assumption $(b_0)$
{is} satisfied. Then
\(
0 < \lambda_* < \overline{\lambda}<\lambda^*
< \lambda_1(\Omega)<\lambda_1(\Omega_0).
\)
For each
\(
\lambda\in(\overline{\lambda},\lambda_1(\Omega_0)),
\)
Problem~\eqref{eq1} admits at least one
{positive solution}.
{
Moreover, there exists $\zeta>0$,    such that
$\overline{\lambda}<\lambda^*-\zeta$ {and, for every} 
$\lambda\in(\lambda^*-\zeta,\lambda_1(\Omega_0))${, Problem~\eqref{eq1}admits}
 a positive minimizing solution
$u_\lambda\in\mathcal{N}_\lambda^+$.}
\noindent For this minimizing solution
 the following properties hold:

\begin{itemize}
  \item[$(i)$]  {$E_{\lambda}(u_\lambda) > 0$, for each $\lambda \in (\lambda^*-\zeta,\lambda^*)$;}
    \item[$(ii)$] $E_{\lambda^*}(u_{\lambda^*}) = 0$;
    \item[$(iii)$] $E_{\lambda}(u_\lambda) < 0$, for each $ \lambda^*<\lambda<\lambda_1(\Omega_0)$;
    \item[$(iv)$] $\|u_\lambda\| \to \infty$, as $\lambda \to \lambda_1(\Omega_0)^-$;
    \item[$(v)$] {the map $\lambda\mapsto c_\lambda^+$ is continuous and strictly decreasing on $(\lambda^*-\zeta,\lambda_1(\Omega_0))$;}
    \item[$(vi)$] { \(
\frac{u_\lambda}{\|u_\lambda\|}
\longrightarrow
\frac{\varphi_1^0}{\|\varphi_1^0\|} ~ {\text{and} ~ \frac{\|u_\lambda\|_2^2}{\|u_\lambda\|^2}
\longrightarrow
\frac{1}{\lambda_1(\Omega_0)}} 
~\text{strongly in }H_0^1(\Omega),
\)
as $\lambda\to\lambda_1(\Omega_0)^-$.}
\end{itemize}
\end{theorem}
\begin{remark}
 Although the condition on $||b||_\infty$ in $(b_0)$ ensures that $\lambda_*>0$, we emphasize that { the key requirement for the validity of the results established in this work is precisely the positivity of $\lambda_*$.} Therefore, alternative assumptions may be considered, provided {that} they  ensure this property. The existence of solutions in the case $\lambda_* \leq 0$ remains an open problem and warrants further investigation.
\end{remark}

{The Rayleigh quotient approach developed above provides, in particular, the minimizing solution \(u_\lambda\in\mathcal N_\lambda^+\) for \(\lambda\in(\lambda^*-\zeta,\lambda_1(\Omega_0))\), thereby yielding the mountain-pass geometry in the range \((\lambda^*-\zeta,\lambda_1(\Omega))\). Combined with Proposition~\ref{low-multiplicity}, this enables us to establish multiplicity throughout the entire interval \((\bar\lambda,\lambda_1(\Omega))\).

\begin{theorem}\label{t1}
Suppose that \(2<q<p<2^*\) and that assumption \((b_0)\) is satisfied. Then, for every
\(
\bar\lambda<\lambda<\lambda_1(\Omega),
\)
Problem~\eqref{eq1} admits at least two distinct positive solutions.

Moreover, if
\(
\lambda^*-\zeta<\lambda<\lambda_1(\Omega),
\)
then one of these solutions is the minimizing solution
\(u_\lambda\in\mathcal N_\lambda^+\) provided by Theorem~\ref{t2}, whereas the other is a positive mountain-pass solution \(v_\lambda\). In particular,
\(
E_\lambda(v_\lambda)=c_\lambda^m>0,
\)
and
\(
c_\lambda^m\to0
~\text{as }\lambda\to\lambda_1(\Omega)^-.
\)
\end{theorem}}
The next theorem provides an existence result for the Problem $(\color{blue}{P_{\overline{\lambda}}})$.
\begin{theorem}\label{t3}
    Suppose that the assumptions $2<q<p<2^*$ and $(b_0)$ are satisfied. If $\lambda=\overline{\lambda}$, then Problem \eqref{eq1} admits at least one {positive solution} $u_\lambda$.
\end{theorem}
Moreover, there exists a range of the parameter $\lambda$ for which the problem admits no positive solution.
\begin{theorem}\label{t4}
Assume $2 < q < p < 2^*$ and that $(b_0)$ holds. If $\lambda < \overline{\lambda}$, then Problem \eqref{eq1} has no nontrivial solution. Moreover, \eqref{eq1} admits no positive solution for any $\lambda\geq \lambda_1(\Omega_0)$.
\end{theorem}
{
\begin{remark}
The previous results show that $\lambda_1(\Omega_0)$ is the sharp upper threshold for the existence of positive solutions. Indeed, Problem~\eqref{eq1} admits {at least one } positive solution for every
\(
\bar\lambda<\lambda<\lambda_1(\Omega_0),
\)
whereas no positive solution exists for $\lambda\geq\lambda_1(\Omega_0)$. Moreover, the minimizing solutions satisfy
\(
\|u_\lambda\|\to\infty
~\text{as }\lambda\to\lambda_1(\Omega_0)^-.
\)
\end{remark}
}


{Based on the previous results, the following schematic bifurcation diagrams
summarize the behavior of the solutions obtained in this work.}

\begin{figure}[!ht]
\centering

\begin{minipage}[t]{0.48\textwidth}
    \centering
    \includegraphics[width=\textwidth]{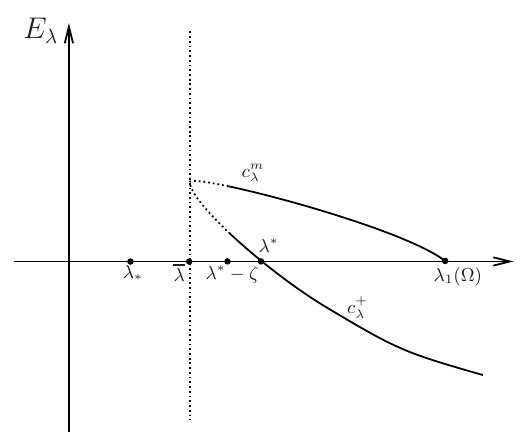}
    \caption{Schematic energy diagram.}
    \label{bif-E}
\end{minipage}
\hfill
\begin{minipage}[t]{0.48\textwidth}
    \centering
    \includegraphics[width=\textwidth]{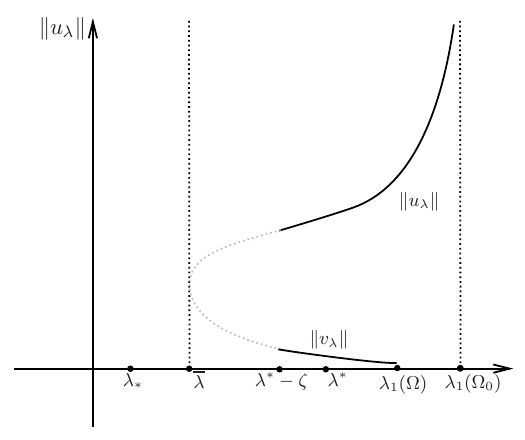}
    \caption{Schematic norm diagram.}
    \label{bif-n}
\end{minipage}

\end{figure}
{Figures \ref{bif-E} and \ref{bif-n} provide schematic summaries of the established parameter ranges, energy signs, and asymptotic behavior. They do not imply the existence of a continuous or connected branch of solutions over the interval $[\bar{\lambda}, \lambda^*-\zeta]$.}

This paper is organized as follows: Section \ref{SB2} introduces the fundamental concepts of the Rayleigh quotient and the Nehari manifold methods, establishing the theoretical foundations for the development of the main theorems and the investigation of the solutions to the Problem \eqref{eq1}. Section \ref{SB3} addresses the proofs of the existence of solutions. Finally, Section \ref{SB4} is devoted to the analysis of the asymptotic behavior of the obtained solutions and we present the proofs of our main results. \\

Throughout this work we shall use the following notations:
	
	\begin{itemize}
		\item The real number $C$ denotes positive constants (possibly different). 
		\item The norm in $L^{p}(\Omega)$ and $L^{\infty}(\Omega)$, will be denoted respectively by $\|\cdot\|_{p}$ and $\|\cdot\|_{\infty}, ~ p \in [1, \infty)$.
		\item We also write $||u||_{r,b}^r:= \Int b(x)|u|^rdx$. 
		\item {$S_r$ denotes the best constant for the embedding $H_0^1(\Omega)\hookrightarrow L^r(\Omega)$, for each $r \in [2, 2^*)$.}
	\end{itemize}

	\section{Preliminary Results}\label{SB2}
 In this section, we present the fundamental concepts of the Rayleigh quotient and the Nehari manifold methods, establishing the theoretical foundations for the development of the main theorems and the investigation of the solutions to the Problem \eqref{eq1}.
\subsection{The Rayleigh Quotient: Definitions and Fundamental Properties}

 In view of {the} identities \eqref{Rn,t} and \eqref{Re,t}, we get
  $$||u||^2_2\dfrac{d}{dt}R_{n}(tu)=(p-2)t^{p-3}||u||^{p}_{p,b}-(q-2)t^{q-3}||u||_{q,b}^q,\  
  \ \ u \in D,\ t>0,$$  {and} $$\frac{1}{2}||u||^2_2\dfrac{d}{dt}R_e(tu)=\dfrac{(p-2)}{p}t^{p-3}||u||_{p,b}^p-\dfrac{(q-2)}{q}t^{q-3}||u||_{q,b}^q
  , \ \ u\in D,\ t>0.$$
 
  It is easy to see that $\dfrac{d}{dt}R_{n}(tu)=0$, if and only if, 
  \begin{equation}\label{tn}
      \tau(u):=t= \left[\dfrac{q-2}{p-2}\dfrac{||u||_{q,b}^q}{||u||^p_{p,b}}\right]^{\frac{1}{p-q}}.
  \end{equation}

  {Similarly}, we get  $\dfrac{d}{dt}R_e(tu)=0$, if and only if, \begin{equation}\label{te}
      {s}(u):=t= \left[\dfrac{p}{q}\right]^{\frac{1}{p-q}}\left[\dfrac{q-2}{p-2}\dfrac{||u||_{q,b}^q}{||u||^p_{p,b}}\right]^{\frac{1}{p-q}}=\left[\dfrac{p}{q}\right]^{\frac{1}{p-q}}\tau(u).
  \end{equation}
  Note that $p>q$ implies $s(u)>\tau(u)$.

\begin{remark} \label{rmk2}
The equation below will be useful for comparing the fibers of the quotients $R_n$ and $R_e$:
\begin{equation*}
R_n(tu)-R_e(tu)=\dfrac{t}{2}\dfrac{d}{dt}R_e(tu), \ \ \text{for each} \ t > 0.
\end{equation*}
  \end{remark}

By using \eqref{Rn,t}, \eqref{Re,t}, \eqref{tn} and \eqref{te}, 
we may conclude that 
	\begin{equation}\label{Qtc}
	  \lambda_n(u):=R_n(\tau(u)u)= \dfrac{||u||^2}{||u||_2^2}-
      {C_{p,q}}\mathcal{R}(u)
	\end{equation} and
 \begin{equation}\label{lambdae}
             \lambda_e(u):=R_e(s(u)u)= \dfrac{||u||^2}{||u||_2^2}-\dfrac{2}{p}\left[\dfrac{p}{q}\right]^\frac{p-2}{p-q}
             {C_{p,q}}\mathcal{R}(u),
	\end{equation}
	with
	\begin{equation}\label{const}
		C_{p,q} := \left(\dfrac{q-2}{p-2}\right)^{\frac{{q-2}}{p-q}}\dfrac{{p-q}}{p-2} \ \text{ and }\  \mathcal{R}(u):=\dfrac{1}{||u||^2_2}\dfrac{||u||_{q,b}^{q\frac{p-2}{p-q}}}{||u||^{p\frac{q-2}{p-q}}_{p,b}}.
	\end{equation}
 
	Furthermore, we observe that
	$\dfrac{d}{dt}R_n(tu)<0$, for each $0<t<\tau(u)$, and $\dfrac{d}{dt}R_n(tu)>0$, for each $ t>\tau(u)$. This can be achieved by utilizing the fact that  $\tau(u)$ is the unique critical point of $t\mapsto R_n(tu)$. Moreover, from \eqref{Rn,t}  we obtain
	\begin{equation}\label{limites}
		\lim_{t\to+\infty}R_n(tu)=+\infty\ \text{and} \  R_n(0)= \dfrac{||u||^2}{||u||_2^2},
	\end{equation}  
  which allow us to state that
$$\lambda_n(u)=\min_{t>0}R_n(tu)=R_n(\tau(u)u).$$
  Using the same argument for \( R_e \), we deduce that 
$$\lambda_e(u)=\min_{t>0}R_e(tu)=R_e(s(u)u).$$

 According to the relationship established in Remark \ref{rmk2}, we deduce that 
 \begin{equation}\label{comp Rn, Re}
 R_n(tu) \leq R_e(tu), ~\text{for all} ~ 0<t\leq s(u) ~\text{and} ~ R_n(tu) > R_e(tu), \text{for all} ~ t>s(u).
 \end{equation}
    
A key step in applying the techniques developed in this work to obtain solutions to Problem \eqref{eq1} is to show that the values $\lambda_*$ and $\lambda^*$, as defined in \eqref{lambdane}, are actually attained within the set $D$. The proofs of these facts are given in Propositions \ref{ln} and \ref{le} below. Before stating them, we establish an important relation between the first eigenvalue of the Laplacian operator and these values.
\begin{lemma}\label{uinD}
    If $2<q<p<2^*$ then $\lambda_*, \lambda^*<\lambda_1(\Omega)$.  \end{lemma}
    \begin{proof} 
    Consider the first eigenfunction $\varphi_1>0$ corresponding to $\lambda_1(\Omega)$, the first eigenvalue of the operator $-\Delta$ in $H_0^1(\Omega)$.  
Since $\varphi_1>0$, we have {$\varphi_1 \in D$}.  
 Moreover, using \eqref{Qtc}, 
the constant $C_{p,q}$ and $\mathcal{R},$ given in \eqref{const}, we obtain 
\[
\lambda_1(\Omega) = \frac{\|\varphi_1\|^2}{\|\varphi_1\|_2^2} > \frac{\|\varphi_1\|^2}{\|\varphi_1\|_2^2}-\dfrac{C_{p,q}}{\|\varphi_1\|^2_2}\dfrac{\|\varphi_1\|_{q,b}^{q\frac{p-2}{p-q}}}{\|\varphi_1\|^{p\frac{q-2}{p-q}}_{p,b}}=\frac{\|\varphi_1\|^2}{\|\varphi_1\|_2^2} - C_{p,q} \mathcal{R}(\varphi_1)\geq \lambda_*.\] 
If we use \eqref{lambdae} in place of \eqref{Qtc} in the inequality above, we obtain $\lambda_1(\Omega)>\lambda^*$.
\end{proof}

\ignore{\begin{lemma}\label{uinD}
    If $2<q<p<2^*$ and {\color{blue}if} $(b_0)$ {\color{blue}holds}. If $\lambda_*=\displaystyle\inf_{v \in D}\lambda_n(v)$ $($or $\lambda^*=\displaystyle\inf_{v \in D}\lambda_e(v)$$)$ is achieved for a function $u \in H_0^1(\Omega)\setminus \{0\}$, then $u \in D$.
    \begin{proof} We prove the result for $\lambda_*$, the proof for $\lambda^*$ is analogous.
    Consider the first eigenfunction $\phi_1>0$ corresponding to $\lambda_1(\Omega)$, the first eigenvalue of the operator $-\Delta$ in $H_0^1(\Omega)$.  
Since $\phi_1 > 0$, we have that $\phi_1 \in D$.  
 Moreover, using \eqref{Qtc}, 
the constant $C_{p,q}$ and $\mathcal{R},$ given in \eqref{const}, we obtain 
\[
\lambda_1(\Omega) = \frac{\|\phi_1\|^2}{\|\phi_1\|_2^2} > \frac{\|\phi_1\|^2}{\|\phi_1\|_2^2}-\dfrac{C_{p,q}}{\|\phi_1\|^2_2}\dfrac{\|\phi_1\|_{q,b}^{q\frac{p-2}{p-q}}}{\|\phi_1\|^{p\frac{q-2}{p-q}}_{p,b}}=\frac{\|\phi_1\|^2}{\|\phi_1\|_2^2} - C_{p,q} \mathcal{R}(\phi_1)\geq \lambda_*,\] 
This proves that $\lambda_1(\Omega) > \lambda_*.$

Suppose, by contradiction, that $\lambda_*$ is realized by some $u \notin D$.  
In this case,
\[
\|u\|_{p,b} = \|u\|_{q,b} = 0
\]
and also
$$
\lambda_* = R_n(u) = \frac{\|u\|^2}{\|u\|_2^2} \geq \lambda_1(\Omega).
$$ {\color{orange} (Se o domínio de $R_n$ é o conjunto $D$, faz sentido escrever $R_n(u)$ para $u \notin D$?)}
This is a contradiction, since we obtained that
$
\lambda_* < \lambda_1(\Omega).
$ Consequently, $u\in D$.
         \end{proof}
\end{lemma}}
In the next result we establish essential properties of the functional $\lambda_n$, given by \eqref{Qtc}. In particular, we prove the existence of a global minimum $\lambda_*$ for $\lambda_n$ and, by applying Lemma~\ref{uinD}, we show that such a minimum is attained in the set $D$.
 \begin{proposition}\label{ln} Assume $2<q<p<2^*$ and that $(b_0)$ holds. Then we obtain the following statements:
		\begin{itemize}
        \item[$(i)$] $u\mapsto \lambda_n(u)=R_n(\tau(u)u)$ is 0-homogeneous; 
			\item[$(ii)$] $u\mapsto \lambda_n(u)$  is  bounded from below. Furthermore, the functional $u\mapsto \lambda_n(u)$ is unbounded from above, that is, there exists a sequence $(w_k) \subset D$ such that $\lambda_n(w_k) \to +\infty$ as $k\to +\infty;$
			\item[$(iii)$] there exists $u \in D$ such that $\displaystyle \lambda_* =\inf_{v \in D}\lambda_n(v)$; 
            \item[$(iv)$]  considering  the function $u$ obtained in item (iii), we have {that} $u_{\lambda_*} := \tau(u)u$ is a weak solution for the problem
            \begin{equation}\label{eq2}
		\left\{
            \begin{array}{rcl}
			-2\Delta u  &=& 2\lambda_* u-b(x)(p|u|^{p-2}u- q|u|^{q-2}u)~\mbox{in}   \  \Omega,  \\[2ex]
			u&=&0 ~ \mbox{on} \ \partial \Omega.
		\end{array}
		\right.
	\end{equation}
     
		\end{itemize}
	\end{proposition}
 \begin{proof} The proof of item (i) follows directly from the definition of $\lambda_n(u)$, given by \eqref{Qtc}. We proceed to prove item (ii): by \eqref{Qtc} {and \eqref{const}} we obtain \begin{equation}\label{a1}\lambda_n(u)= \dfrac{\|u\|^2}{\|u\|_2^2}-
 {C_{p,q}}\mathcal{R}(u)
 \geq \lambda_1(\Omega)-C_{p,q}\mathcal{R}(u). 
 \end{equation}   
 Since $2<q<p$, by using the interpolation law, we write
\begin{equation}\label{interpol}
\|u\|_{q,b}\leq \|u\|_{2,b}^{\frac{2}{q}\frac{p-q}{p-2}}\|u\|_{p,b}^{\frac{p}{q}\frac{q-2}{p-2}}\leq \|b\|_{_\infty}^{{\frac{1}{q}\frac{p-q}{p-2}}}\|u\|_2^{\frac{2}{q}\frac{p-q}{p-2}}\|u\|_{p,b}^{\frac{p}{q}\frac{q-2}{p-2}}.
\end{equation}  It follows from \eqref{const} and \eqref{interpol} that 
$$\mathcal{R}(u)\leq \|b\|_\infty.$$
{
From the first inequality in \eqref{interpol} and the definition of $\mathcal{R}$ in \eqref{const},
\begin{equation}\label{R-weighted-L2}
\mathcal{R}(u)
\leq \frac{\|u\|_{2,b}^{2}}{\|u\|_2^{2}},\qquad u\in D.
\end{equation}
}
As a consequence  of the above inequality, \eqref{a1} and hypothesis $(b_0)$, we may conclude that
{\begin{equation}\label{lambda*positivo}\lambda_n(u)\geq \lambda_1(\Omega)-C_{p,q}\|b\|_\infty>0.\end{equation}}

It follows immediately that, $\lambda_*=\displaystyle\inf_{v\in D}\lambda_n(v)>0.$

 {In order to verify that $u\mapsto \lambda_n(u)$ is unbounded from above, choose a ball
$B\subset\subset\Omega\setminus\overline{\Omega_0}$, which exists because
$\Omega_0\subset\subset\Omega$. By assumption $(b_0)$, $b>0$ a.e. in $B$. Let $(w_k)$ be an
orthonormal sequence in $H_0^1(B)$ with respect to the Dirichlet inner product,
and extend each $w_k$ by zero outside $B$. Then $w_k\in D$, $\|w_k\|=1$, and
$w_k\rightharpoonup0$ in $H_0^1(\Omega)$. By the compact embedding
$H_0^1(\Omega)\hookrightarrow L^2(\Omega)$, we have $\|w_k\|_2\to0$. Moreover, since
$\mathcal R(w_k)\leq\|b\|_\infty$, it follows that
\[
\lambda_n(w_k)
=\frac{1}{\|w_k\|_2^2}-C_{p,q}\mathcal R(w_k)
\geq\frac{1}{\|w_k\|_2^2}-C_{p,q}\|b\|_\infty
\longrightarrow+\infty.
\]
Thus $\lambda_n$ is unbounded from above.}

Now, we shall prove (iii): 
  Let $(u_k) \subset D, \|u_k\|=1$ be a minimizing sequence of $\lambda_n$ 
  i.e.,
 $   \lambda_n(u_k) = \lambda_*+o_k(1).$
Since $(u_k)$ is bounded, there exists $u \in H_0^1(\Omega)$ such that $u_k \rightharpoonup u$ weakly. As a consequence, using the compact embedding $H^1_0(\Omega) \hookrightarrow  L^r(\Omega)$ for $r \in [2,2^*)$,  the fact that $D \subset H^1_0(\Omega)$ and the {sequential lower semicontinuity of the norm}, we obtain
\begin{equation}\label{Se}
    \|u\| \leq \displaystyle\liminf_{{k\to\infty}} \|u_k\|, \ \ 
    \displaystyle\lim_{k\to\infty}\|u_k\|_{p,b} =\|u\|_{p,b}\ \text{and} \  \displaystyle\lim_{k\to\infty}\|u_k\|_{q,b} =\|u\|_{q,b}.
\end{equation}
\textbf{Claim:} $u \neq 0.$
Suppose, by contradiction, that $u = 0$. Then we have $u_k \rightharpoonup 0$ weakly in $H_0^1(\Omega)$. It follows from \eqref{a1} and preceding arguments  that
\begin{equation*}
\lambda_* =  \displaystyle\liminf_{{k\to\infty}}\lambda_n(u_k)
\geq \displaystyle\liminf_{{k\to\infty}}\dfrac{1}{\|u_k\|_2^2} - C_{p,q} {\| b \|_{\infty}}. 
\end{equation*}
{
Since $u_k\rightharpoonup 0$ in $H_0^1(\Omega)$ and the embedding
$H_0^1(\Omega)\hookrightarrow L^2(\Omega)$ is compact, we have
$u_k\to0$ strongly in $L^2(\Omega)$. Hence $\|u_k\|_2\to0$, and the preceding lower bound implies
$\lambda_n(u_k)\to+\infty$, contradicting
$\lambda_n(u_k)=\lambda_*+o_k(1)$.
Thus, $u\neq0$, which verifies the claim.
} 

From \eqref{a1} and \eqref{Se}, we now observe that
\begin{eqnarray*}
    \lambda_* = \inf_{v \in D } \lambda_n(v) &=&  \displaystyle\liminf_{{k\to\infty}} \lambda_n(u_k)
=\displaystyle\liminf_{{k\to\infty}}\left[\dfrac{\| u_k \|^2}{\|u_k\|^2_2} - C_{p,q} \mathcal{R}(u_k)\right]\\&\geq &\left[\dfrac{\| u \|^2}{\|u\|^2_2} - C_{p,q} \mathcal{R}(u)\right]=\lambda_n(u)\geq\lambda_*, ~\text{if} \ u \in D.
\end{eqnarray*} 
 {
On the other hand, suppose that $u\notin D$. By the definition of $D$ and hypothesis $(b_0)$, we have 
$u=0$ a.e. in $\Omega\setminus\Omega_0$, and therefore
$\|u\|_{2,b}=0$. Since $u_k\to u$ strongly in $L^2(\Omega)$ and $b\in L^\infty(\Omega)$, it follows that
$
\|u_k\|_{2,b}\longrightarrow \|u\|_{2,b}=0.
$
Moreover, the preceding claim yields $u\neq0$, so that $\|u_k\|_2\to\|u\|_2>0$. By \eqref{R-weighted-L2}, we obtain
\[
0\leq \mathcal R(u_k)
\leq \frac{\|u_k\|_{2,b}^2}{\|u_k\|_2^2}
\longrightarrow0.
\]
Consequently,
\[
\lambda_*
=\liminf_{k\to\infty}\lambda_n(u_k)
\geq \frac{\|u\|^2}{\|u\|_2^2}
\geq \lambda_1(\Omega),
\]
contradicting Lemma~\ref{uinD}, which establishes that $\lambda_*<\lambda_1(\Omega)$. Therefore $u\in D$, and the infimum $\lambda_*$ is attained at $u\in D$.
}

For the proof of $(iv)$, let $v:=u_{\lambda_*}=\tau(u)u$. Since $\lambda_n$ is $0$-homogeneous and $u$ minimizes $\lambda_n$, the function $v$ also minimizes $\lambda_n$, and $\tau(v)=1$. Moreover,
\[
\frac{d}{dt}R_n(tv)\Big|_{t=1}=R_n'(v)[v]=0.
\]
{Since $\tau$ is $C^1$ on $D$, the chain rule applied to
$\lambda_n(z)=R_n(\tau(z)z)$ yields, for every $\varphi\in H_0^1(\Omega)$,
\[
0=\lambda_n'(v)[\varphi]
 =R_n'(v)\big[\tau'(v)[\varphi]v+\varphi\big]
 =R_n'(v)[\varphi],
\]
where the last equality follows from $R_n'(v)[v]=0$. Thus $v$ is a critical
point of the quotient $R_n=\widehat F/\widehat G$, where
\[
\widehat F(z)=\|z\|^2+\|z\|_{p,b}^p-\|z\|_{q,b}^q,
\qquad
\widehat G(z)=\|z\|_2^2.
\]
Since $R_n(v)=\lambda_*$, differentiation of the quotient gives
\[
0=R_n'(v)[\varphi]
=\frac{1}{\widehat G(v)}
\big(\widehat F'(v)[\varphi]-\lambda_*\widehat G'(v)[\varphi]\big).
\]
Consequently,
\[
-2\Delta v
=2\lambda_*v-b(x)\big(p|v|^{p-2}v-q|v|^{q-2}v\big)
\]
in the weak sense, which is precisely Problem~\eqref{eq2}.}
\end{proof}

The following result establishes conclusions analogous to those of Proposition \ref{ln}, but with respect to the functional $\lambda_e$, given by \eqref{lambdae}. Moreover, it ensures the existence of a solution to Problem \eqref{eq1} in the particular case $\lambda = \lambda^*$, namely a solution with zero energy.


\begin{proposition}\label{le} Assume $2<q<p<2^*$ and that $(b_0)$ holds. Then we obtain the following statements:
		\begin{itemize}
        \item[$(i)$] $u\mapsto \lambda_e(u)=R_e(s(u)u)$ is 0-homogeneous; 
	\item[$(ii)$] $u\mapsto R_e(s(u)u)$  is  bounded from below;
			\item[$(iii)$] there exists $u \in D$ such that $\displaystyle \lambda^* :=\inf_{v\in D}\inf_{t>0}R_e(tv):= \inf_{v \in D} \lambda_e(v)$;
            
            \item[$(iv)$] considering the function $u$ obtained in item $(iii)$ we have {that} $u_{\lambda^*} := s(u)u
            $ 
            is a weak solution for the problem
            \begin{equation}\label{eq3}
		\left\{
            \begin{array}{rcl}
			-\Delta u  &=& \lambda^* u-b(x)(|u|^{p-2}u- |u|^{q-2}u)~\mbox{in}   \  \Omega,  \\[2ex]
			u&=&0 ~ \mbox{on} \ \partial \Omega.
		\end{array}
		\right.
	\end{equation}
    \end{itemize}
  \end{proposition}
  \begin{proof}
The proofs of {$(i)$--$(ii)$} follow the same ideas as
those of Proposition~\ref{ln}, and hence we omit the details.
{
To prove $(iii)$, set
\[
C^e_{p,q}:=\frac{2}{p}\left(\frac{p}{q}\right)^{\frac{p-2}{p-q}}C_{p,q}.
\]
Then, by \eqref{lambdae},
\[
\lambda_e(v)=\frac{\|v\|^2}{\|v\|_2^2}-C^e_{p,q}\mathcal R(v),\qquad v\in D.
\]
Let $(u_k)\subset D$, with $\|u_k\|=1$, be a minimizing sequence for $\lambda_e$. Passing to a subsequence, we may assume that
$u_k\rightharpoonup u$ in $H_0^1(\Omega)$ and $u_k\to u$ in $L^r(\Omega)$ for every $2\leq r<2^*$. If $u=0$, then $\|u_k\|_2\to0$. Since $\mathcal R(u_k)\leq\|b\|_\infty$, we obtain 
\[
\lambda_e(u_k)\geq\frac1{\|u_k\|_2^2}-C^e_{p,q}\|b\|_\infty\longrightarrow+\infty,
\]
which contradicts the fact that $(u_k)$ is a minimizing
sequence. Thus, $u\neq0$.

Suppose now that $u\notin D$. Then $\|u\|_{2,b}=0$. Since $u_k\to u$ strongly in $L^2(\Omega)$, we have $\|u_k\|_{2,b}\to0$ and $\|u_k\|_2\to\|u\|_2>0$. Consequently, \eqref{R-weighted-L2} yields $\mathcal R(u_k)\to0$. It follows that
\[
\lambda^*
=\liminf_{k\to\infty}\lambda_e(u_k)
\geq\frac{\|u\|^2}{\|u\|_2^2}
\geq\lambda_1(\Omega),
\]
contradicting Lemma~\ref{uinD}, which asserts that $\lambda^*<\lambda_1(\Omega)$. Therefore $u\in D$. By the strong convergence in the weighted $p$- and $q$-norms
and in $L^2(\Omega)$, together with the weak lower
semicontinuity of the Dirichlet norm, we obtain $
\lambda^*
\geq\lambda_e(u).
$
On the other hand, since $u\in D$ and $\lambda^*$ is the
infimum of $\lambda_e$ over $D$, we have
$\lambda_e(u)\geq\lambda^*$. Hence $
\lambda^*\geq\lambda_e(u)\geq\lambda^*,
$
and the infimum in $(iii)$ is attained.
}

To prove $(iv)$, let $v:=u_{\lambda^*}=s(u)u$. Since $\lambda_e$ is $0$-homogeneous and $u$ minimizes $\lambda_e$, the function $v$ also minimizes $\lambda_e$, and $s(v)=1$. Moreover,
\[
\frac{d}{dt}R_e(tv)\Big|_{t=1}=R_e'(v)[v]=0.
\]
{Since $s$ is $C^1$ on $D$, the chain rule applied to
$\lambda_e(z)=R_e(s(z)z)$ provides, for every $\varphi\in H_0^1(\Omega)$,
\[
0=\lambda_e'(v)[\varphi]
=R_e'(v)\big[s'(v)[\varphi]v+\varphi\big]
=R_e'(v)[\varphi].
\]
Writing
\[
F(z)=\frac12\|z\|^2+\frac1p\|z\|_{p,b}^p-\frac1q\|z\|_{q,b}^q,
\qquad
G(z)=\frac12\|z\|_2^2,
\]
we note that $R_e=F/G$ and $R_e(v)=\lambda^*$. Hence
\[
0=R_e'(v)[\varphi]
=\frac1{G(v)}\big(F'(v)[\varphi]-\lambda^*G'(v)[\varphi]\big)
=\frac1{G(v)}E'_{\lambda^*}(v)[\varphi].
\]
Therefore $E'_{\lambda^*}(v)=0$, which implies that $v=u_{\lambda^*}$ is a weak solution of Problem~\eqref{eq3}.}
 \end{proof}
By Propositions \ref{ln} and \ref{le}, it follows that \eqref{lambdane} can be expressed in the form 
$$\lambda_*=\inf_{u\in D}\lambda_n(u)\qquad\mbox{and}\qquad \lambda^*=\inf_{u\in D}\lambda_e(u).$$
\begin{cor}
   Suppose $2<q<p<2^*$ and that $(b_0)$ holds. Then $ \lambda_*<\lambda^*$.
   \begin{proof}
        It follows from item $(iii)$ of Propositions \ref{ln} and \ref{le} that there exist $u_{*}, u^*\in D$ such that $\lambda_*=R_n(\tau(u_*)u_*)$ and $\lambda^*=R_e(s(u^*)u^*)$.  Combining these relations with the inequality  $\tau(u)<s(u)$, the monotonicity of $R_n$ for $t>\tau(u)$ and \eqref{comp Rn, Re}, we deduce that
       $$\lambda_*\leq R_n(\tau(u^*)u^*)<R_n(s(u^*)u^*)\leq R_e(s(u^*)u^*)=\lambda^*.$$
   \end{proof}
\end{cor}
The proofs of the next two propositions are straightforward and follow {the} arguments from \cite{yavdat0,yavdat1,yavdat2}.
\begin{proposition}\label{1}  
Let \( t > 0 \) and \( u \in D \) be fixed. Then, the condition \( tu \in \mathcal{N}_{\lambda} \) is equivalent to \( R_n(tu) = \lambda \).  
More specifically, {we have} \( R_n(tu) = \lambda \) if and only if \( E'_{\lambda}(tu)tu = 0 \). In this case, we have
\begin{equation}\label{import1}
       \dfrac{d}{dt}R_n(tu)=\dfrac{1}{t}\dfrac{E''_\lambda(tu)(tu)^2}{||tu||_2^2}.\end{equation}
Additionally, the following properties hold:  
\begin{itemize}  
    \item[$(i)$] \( R_n(tu) > \lambda \) if and only if \( E'_{\lambda}(tu)tu > 0 \);  
    \item[$(ii)$] \( R_n(tu) < \lambda \) if and only if \( E'_{\lambda}(tu)tu < 0 \).  
\end{itemize}  
\end{proposition}  

\begin{proposition}\label{2}  
Let \( u \in D \) and \( t > 0 \) be fixed. Then, {we have} \( R_e(tu) = \lambda \) if and only if \( E_{\lambda}(tu) = 0 \).  In this case, we have
 \begin{equation}\label{import2}
       \dfrac{d}{dt}R_e(tu)=\dfrac{2}{t}\dfrac{E'_\lambda(tu)(tu)}{||tu||_2^2}.\end{equation}
Moreover, the following properties hold:  
\begin{itemize}  
    \item[$(i)$] \( R_e(tu) > \lambda \) if and only if \( E_{\lambda}(tu) > 0 \);  
    \item[$(ii)$] \( R_e(tu) < \lambda \) if and only if \( E_{\lambda}(tu) < 0 \).  
     
\end{itemize}  
\end{proposition}  

\begin{figure}[!ht]
\centering
\IfFileExists{fibra_R.pdf}{%
    \includegraphics[scale=0.7]{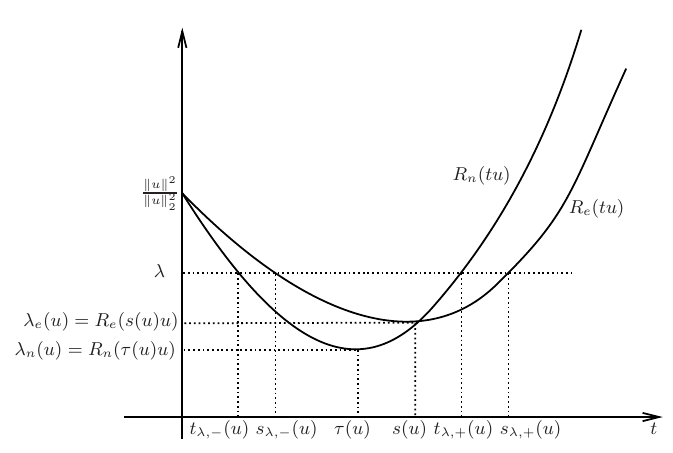}%
}{%
\begin{tikzpicture}[x=4.0cm,y=1.15cm]
    \draw[->] (0,0) -- (1.62,0) node[right] {$t$};
    \draw[->] (0,0) -- (0,4.95) node[above] {$R(tu)$};
    \draw[thick,domain=0:1.47,samples=120,smooth]
        plot (\x,{1 + (\x)^4 - (\x)});
    \draw[thick,dashed,domain=0:1.47,samples=120,smooth]
        plot (\x,{1 + (1/3)*(\x)^4 - (2/3)*(\x)});
    \node[anchor=west] at (1.17,2.22) {$R_n(tu)$};
    \node[anchor=west] at (1.18,1.04) {$R_e(tu)$};
    \draw[densely dotted] ({pow(1/4,1/3)},0) -- ({pow(1/4,1/3)},{1 + pow(pow(1/4,1/3),4)-pow(1/4,1/3)});
    \draw[densely dotted] ({pow(1/2,1/3)},0) -- ({pow(1/2,1/3)},{1 + (1/3)*pow(pow(1/2,1/3),4)-(2/3)*pow(1/2,1/3)});
    \node[anchor=north east] at ({pow(1/4,1/3)},-0.03) {$\tau(u)$};
    \node[anchor=north west] at ({pow(1/2,1/3)},-0.03) {$s(u)$};
    \fill ({pow(1/4,1/3)},{1 + pow(pow(1/4,1/3),4)-pow(1/4,1/3)}) circle (0.7pt);
    \fill ({pow(1/2,1/3)},{1 + (1/3)*pow(pow(1/2,1/3),4)-(2/3)*pow(1/2,1/3)}) circle (0.7pt);
\end{tikzpicture}
}
\caption{Schematic behavior of the fibering functions $R_n(tu)$ and $R_e(tu)$. In particular, $\tau(u)<s(u)$, and the two fibers meet again at $t=s(u)$.}
\label{fig1}
\end{figure}

{The following proposition establishes  conditions for a given function $u \in D$ to admit a projection onto the Nehari manifold $\mathcal{N}_\lambda$.}
\begin{proposition}\label{projecao}
Suppose $2<q<p<2^*$ and that $(b_0)$ holds. Given  $u \in D$, we obtain:
\begin{enumerate}
\item[$(1)$] If $\lambda_n(u)<\lambda<\dfrac{||u||^2}{||u||^2_2}$, then the fibering map $t\mapsto E_{\lambda}(tu)$ has exactly two distinct critical points $0 < t_{\lambda,-}(u) < \tau(u) < t_{\lambda,+}(u)$. In addition, the following assertions are valid:

\begin{itemize}
    \item[$(i)$] {The point} $t_{\lambda,+}(u)$ is a local minimum for the fibering map $t\mapsto E_{\lambda}(tu)$ which satisfies $t_{\lambda,+}(u)u \in\mathcal{N}^+_{\lambda}\cap D$. Furthermore, {the point} $t_{\lambda,-}(u)$ is a local maximum for the fibering map $t\mapsto E_{\lambda}(tu)$ which verifies $t_{\lambda,-}(u)u \in \mathcal{N}^-_{\lambda}\cap D$.
    {\item[$(ii)$] Define
    \[
    U_{\lambda,+}:=\{u\in D:\lambda_n(u)<\lambda\},\qquad
    U_{\lambda,-}:=\left\{u\in D:\lambda_n(u)<\lambda<\frac{\|u\|^2}{\|u\|_2^2}\right\}.
    \]
    Then $u\mapsto t_{\lambda,+}(u)$ belongs to $C^1(U_{\lambda,+},\mathbb R)$, whereas $u\mapsto t_{\lambda,-}(u)$ belongs to $C^1(U_{\lambda,-},\mathbb R)$.
    \item[$(iii)$] For fixed $u\in D$, the map $\lambda\mapsto t_{\lambda,+}(u)$ is $C^1$ on $(\lambda_n(u),\infty)$, while $\lambda\mapsto t_{\lambda,-}(u)$ is $C^1$ on
    $\left(\lambda_n(u),\|u\|^2/\|u\|_2^2\right)$. On their respective domains,
    \[
    \frac{d}{d\lambda}t_{\lambda,\pm}(u)=
    \frac{t_{\lambda,\pm}(u)\,\|t_{\lambda,\pm}(u)u\|_2^2}
    {E''_\lambda(t_{\lambda,\pm}(u)u)(t_{\lambda,\pm}(u)u)^2}.
    \]}
    \end{itemize}
    {\item[$(2)$] If $\lambda < \lambda_n(u)$, the fibering map $t\mapsto E_{\lambda}(tu)$ does not admit any critical point.
    \item[$(3)$] If $\lambda = \lambda_n(u) $, the fibering map $t\mapsto  E_{\lambda}(tu)$ has a unique critical point $\tau(u) > 0$ which satisfies $\tau(u)u \in \mathcal{N}^0_{\lambda_n(u)}\cap D$.
    \item[$(4)$] In addition, if $\dfrac{||u||^2}{||u||_2^2}\leq\lambda <\infty$, then the fibering map $t\mapsto E_\lambda(tu)$ has a unique critical point $t_{\lambda,+}(u) > \tau(u)$ which verifies $t_{\lambda,+}(u)u \in\mathcal{N}^+_{\lambda}\cap D$}.
    \end{enumerate}
\end{proposition}
\begin{proof}
(1) Let $u \in D$ be fixed and assume that $R_n(\tau(u)u) < \lambda< \frac{||u||^2}{||u||^2_2}$. Observe that {the function} $t \mapsto R_n(tu), t > 0,$ is  continuous. Moreover, we have shown that $R_n(tu)$ is strictly decreasing for $t \in (0, \tau(u))$ and strictly increasing for $t \in (\tau(u), \infty)$. Thus, considering the above arguments and \eqref{limites}, we obtain that the equation $R_n(tu) = \lambda$ admits exactly two roots (see Figure \ref{fig1}), which we denote by $t_{\lambda,+}(u)$ and $t_{\lambda,-}(u)$, where $0 < t_{\lambda,-}(u) < \tau(u) < t_{\lambda,+}(u)$. We know that $R_n(tu) = \lambda$ {if and only if} $E'_\lambda(tu)tu = 0$. As a consequence, we deduce that $t_{\lambda,+}(u)$ and $t_{\lambda,-}(u)$ are critical points of {the} map $t\mapsto E_\lambda(tu)$. On the other hand, we know that $\dfrac{d}{dt}R_n(t_{\lambda,+}(u)u) > 0$ and $\dfrac{d}{dt}R_n(t_{\lambda,-}(u)u) < 0$. Under these conditions, using \eqref{import1} we infer that
$$
E''_\lambda(t_{\lambda,+}(u)u)(t_{\lambda,+}(u)u)^2 > 0 \quad \text{and} \quad E''_\lambda(t_{\lambda,-}(u)u)(t_{\lambda,-}(u)u)^2 < 0.
$$
The above inequalities imply that $t_{\lambda,-}(u)u \in \mathcal{N}_\lambda^- \cap D$ and $t_{\lambda,+}(u)u \in \mathcal{N}^+_\lambda \cap D$.\\
To prove $(ii)$, we consider the auxiliary function $\Phi^\pm : (0, +\infty) \times D \rightarrow \mathbb{R}$ given by $\Phi^\pm(t, u) = E_\lambda'(tu)tu$. It follows that $\Phi^\pm(t, u) = 0$ if and only if $tu \in \mathcal{N}_\lambda$.  Furthermore, we observe that $\dfrac{\partial \Phi^\pm}{\partial t} (t, u) \neq 0$ for every 
{$(t,u)\in (0,+\infty)\times D$ such that, $ tu\in \mathcal{N}_\lambda^\pm$. } {Thus, the Implicit Function Theorem \cite{drabek} applies separately to the two nondegenerate branches, yielding the $C^1$ regularity stated in item $(ii)$. Item $(iii)$ follows by differentiating the identity $E_\lambda'(t_{\lambda,\pm}(u)u)t_{\lambda,\pm}(u)u=0$ with respect to $\lambda$ on the corresponding domain.}

(2)  Assume that $\lambda <\lambda_n(u)$, see Figure \ref{fig1}. Hence, using Proposition \ref{1} we get {that the equation} $ R_n(tu)= \lambda$ does not admit any root, since $\lambda<R_n(\tau(u)u)=\min_{t>0}R_n(tu)$. This statement is equivalent to saying that $ E'_{\lambda}(tu)(tu) \neq 0$ for all $t > 0$ and for each $u \in D$. As a consequence, $tu \notin \mathcal{N}_{\lambda}$ for all $t > 0$. This finishes the proof of item (2).

(3) Here we assume that $\lambda = \lambda_n(u)$, see Figure \ref{fig1}. In this case, for each $u \in D$, we obtain that $\dfrac{d}{dt}R_n(\tau(u)u)=0$ (see Proposition \ref{1}) where $\tau(u) > 0$ is the unique critical point of $t\mapsto R_n(tu)$. In view of \eqref{import1}  we conclude that $E''_{\lambda}(\tau(u)u)(\tau(u)u)^2 = 0$. Therefore, we obtain $\tau(u)u \in \mathcal{N}^0_{\lambda_n(u)}\cap D$.

{(4)} Let $u \in D$ be fixed and assume that $\dfrac{||u||^2}{||u||_2^2}\leq\lambda <\infty$.  It follows from the continuity of the function $t \mapsto R_n(tu), t > 0$, from its monotonicity for $t>\tau(u)$ and from \eqref{limites} that the equation $R_n(tu) = \lambda$ admits exactly one root (see Figure \ref{fig1}) $t_{\lambda,+}(u)>\tau(u) $ (see Figure \ref{fig1}). We know that $R_n(tu) = \lambda$ {if and only if} $E'_\lambda(tu)tu = 0$. As a consequence, we deduce that $t_{\lambda,+}(u)$ is a {critical point} of the map $t\mapsto E_\lambda(tu)$. Since $\dfrac{d}{dt}R_n(t_{\lambda,+}(u){u}) > 0$, it follows from  \eqref{import1} that
$E''_\lambda(t_{\lambda,+}(u)u)(t_{\lambda,+}(u)u)^2 > 0$. Consequently $t_{\lambda,+}(u)u\in \mathcal{N}_\lambda^+.$
 \end{proof}

 By employing an argument analogous to that used in Proposition \ref{projecao}, we establish the following result:
\begin{proposition}\label{projecaoE}
Suppose $2<q<p<2^*$ and that $(b_0)$ holds. Given  $u \in D$, we obtain:
\begin{enumerate}
{\item[$(1)$] If $\lambda_e(u)<\lambda<\dfrac{||u||^2}{||u||^2_2}$, then the fibering map $t\mapsto E_{\lambda}(tu)$ has exactly two distinct positive zeros
$0<s_{\lambda,-}(u)<s(u)<s_{\lambda,+}(u)$.
Define
\[
W_{\lambda,+}:=\{u\in D:\lambda_e(u)<\lambda\},\qquad
W_{\lambda,-}:=\left\{u\in D:\lambda_e(u)<\lambda<\frac{\|u\|^2}{\|u\|_2^2}\right\}.
\]
Then $u\mapsto s_{\lambda,+}(u)$ belongs to $C^1(W_{\lambda,+},\mathbb R)$ and
$u\mapsto s_{\lambda,-}(u)$ belongs to $C^1(W_{\lambda,-},\mathbb R)$.}
    \item[$(2)$] If $\lambda < \lambda_e(u)$, the fibering map $t\mapsto E_{\lambda}(tu)$ has no root.
    \item[$(3)$] If $\lambda = \lambda_e(u) $, the fibering map $t\mapsto  E_{\lambda}(tu)$ has a unique root $s(u) > 0$.
   {\item[$(4)$] In addition, if $\dfrac{||u||^2}{||u||_2^2}\leq\lambda <\infty$, then the fibering map $t\mapsto E_\lambda(tu)$ has a unique positive zero $s_{\lambda,+}(u)>s(u)$, i.e. $E_\lambda(s_{\lambda,+}(u)u)=0$.} 
    \end{enumerate}
\end{proposition}
\begin{proof}
    For the proof, it suffices to consider \eqref{import2} and argue as in {Proposition} \ref{projecao}.
\end{proof}

{The behavior of the fibers 
$R_n(tu)$ and $R_e(tu)$, as characterized by Propositions \ref{projecao} and \ref{projecaoE}, is illustrated schematically in Figure \ref{fig1}.}

 As a consequence of Proposition \ref{projecaoE} combined with Proposition \ref{projecao}, we can estimate the sign of the energy functional in functions belonging to the Nehari manifold.
\begin{cor}\label{energ-signal}
Suppose $2<q<p<2^*$ and that $(b_0)$ holds.  Given  $u \in D$, we obtain: 
\begin{enumerate}
    \item[$(1)$] if $\lambda_n(u)<\lambda<\lambda_e(u)$, then $0<E_\lambda(t_{\lambda,+}(u)u)<E_\lambda(t_{\lambda,-}(u)u)$;
    \item[$(2)$] if $\lambda=\lambda_e(u)$, then $0=E_\lambda(t_{\lambda,+}(u)u)<E_\lambda(t_{\lambda,-}(u)u)$;
    \item[$(3)$] if $\lambda_e(u)<\lambda <\frac{\|u\|^2}{\|u\|_2^2}$, then $E_\lambda(t_{\lambda,+}(u)u)<0<E_\lambda(t_{\lambda,-}(u)u)$;
    \item[$(4)$] if $\lambda \geq\frac{\|u\|^2}{\|u\|_2^2}$, then there is no projection in $\mathcal{N}_\lambda^-$ and $E_\lambda(t_{\lambda,+}(u)u)<0$.
\end{enumerate}
\end{cor}
\begin{proof}
    {(1) In this case, it follows from item (1) of Proposition \ref{projecao} that the fibering map $t\mapsto E_{\lambda}(tu)$ admits exactly two distinct critical points $t_{\lambda,-}(u)$ and $t_{\lambda,+}(u)$ which correspond to a local maximum and a local minimum of $E_{\lambda}(tu)$, respectively. Moreover, from Proposition \ref{2}, $\lambda<R_e(s(u)u)<R_e(t_{\lambda,+}(u)u)$ implies that $E_\lambda(t_{\lambda,+}(u)u)>0$. Therefore, $0<E_\lambda(t_{\lambda,+}(u)u)<E_\lambda(t_{\lambda,-}(u)u)$.}
    
    {(2) If $\lambda=\lambda_e(u)$, it follows from item (3) of Proposition \ref{projecaoE} that there exists a unique point $s(u)>0$ that verifies $E_\lambda(s(u)u)=0$. Note that, in this case we have $s(u)=t_{\lambda,+}(u)$, which, by using item (1) of Proposition \ref{projecao} once again, shows that $0=E_\lambda(t_{\lambda,+}(u)u)<E_\lambda(t_{\lambda,-}(u)u)$.}

    {(3) Here, item (1) of Proposition \ref{projecaoE} asserts that  $t\mapsto E_{\lambda}(tu)$ has exactly two distinct  points $0 < s_{\lambda,-}(u) < s(u) < s_{\lambda,+}(u)$, such that $E_\lambda(s_{\lambda,-}(u)u)=E_\lambda(s_{\lambda,+}(u)u)=0.$  Moreover, we already know that $t_{\lambda,-}(u)$ is a local maximum and $t_{\lambda,+}(u)$ is a local minimum of $E_{\lambda}(tu)$, respectively, with $E_\lambda(0)=0$. Hence, $E_\lambda(t_{\lambda,+}(u)u)<0<E_\lambda(t_{\lambda,-}(u)u)$.}

{(4) It follows from item (4) of Proposition \ref{projecaoE} that the fibering $t\mapsto E_\lambda(tu)$ has a unique point $s_{\lambda,+}(u) > s(u)$ which verifies $E_\lambda(s_{\lambda,+}(u)u)=0.$ Moreover, item (4) of Proposition \ref{projecao} establishes {that} the fibering $t\mapsto E_\lambda(tu)$ has a unique critical point $t_{\lambda,+}(u) > \tau(u)$ such that $t_{\lambda,+}(u)u \in\mathcal{N}^+_{\lambda}\cap D$, i.e., $E_\lambda^{'}(t_{\lambda,+}(u)u)=0$ and $E_\lambda^{''}(t_{\lambda,+}(u)u)>0$. Consequently, $E_\lambda(t_{\lambda,+}(u)u)<0$. }
    \end{proof}

\begin{figure}[!ht]
\centering
\IfFileExists{fibras_Ea.pdf}{%
    \includegraphics[scale=0.7]{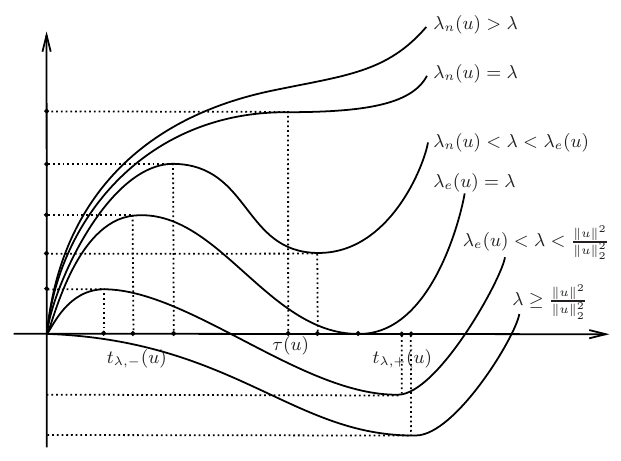}%
}{%
\begin{tikzpicture}[x=3.15cm,y=6.0cm]
    \draw[->] (0,-0.19) -- (1.29,-0.19) node[right] {$t$};
    \draw[->] (0,-0.23) -- (0,0.56) node[above] {$E_\lambda(tu)$};
    \draw[densely dotted] (0,0) -- (1.24,0);

    \draw[thick,domain=0:1.22,samples=120,smooth]
       plot (\x,{0.21*(\x)^2 + (1/6)*(\x)^6 - (1/3)*(\x)^3});
    \draw[thick,dashed,domain=0:1.22,samples=120,smooth]
       plot (\x,{0.198425*(\x)^2 + (1/6)*(\x)^6 - (1/3)*(\x)^3});
    \draw[thick,dotted,domain=0:1.22,samples=120,smooth]
       plot (\x,{0.15*(\x)^2 + (1/6)*(\x)^6 - (1/3)*(\x)^3});
    \draw[thick,dash dot,domain=0:1.22,samples=120,smooth]
       plot (\x,{-0.03*(\x)^2 + (1/6)*(\x)^6 - (1/3)*(\x)^3});

    \draw[thick] (1.38,0.45) -- (1.48,0.45);
    \node[anchor=west] at (1.51,0.45) {$\lambda_n(u)<\lambda<\lambda_e(u)$};
    \draw[thick,dashed] (1.38,0.35) -- (1.48,0.35);
    \node[anchor=west] at (1.51,0.35) {$\lambda=\lambda_e(u)$};
    \draw[thick,dotted] (1.38,0.25) -- (1.48,0.25);
    \node[anchor=west] at (1.51,0.25) {$\lambda_e(u)<\lambda<\dfrac{\|u\|^2}{\|u\|_2^2}$};
    \draw[thick,dash dot] (1.38,0.15) -- (1.48,0.15);
    \node[anchor=west] at (1.51,0.15) {$\lambda\geq\dfrac{\|u\|^2}{\|u\|_2^2}$};
\end{tikzpicture}
}
\caption{Schematic energy fibers associated with the four parameter regimes in Corollary \ref{energ-signal}.}
\label{fig2}
\end{figure}

{The fibering map $t \mapsto E_\lambda(tu)$, characterized by Corollary \ref{energ-signal}, is illustrated schematically in Figure \ref{fig2}.}

\subsection{Structure and Fundamental Properties of $\mathcal{N}_\lambda$}
Here, we present the main characteristics with respect to \(\mathcal{N}_\lambda^-\), \(\mathcal{N}_\lambda^+\) and \(\mathcal{N}_\lambda^0\), which will support the study of the solutions to the Problem \eqref{eq1}.
 We start this section by establishing that Problem \eqref{eq1} admits no solution within a certain parameter range.
\begin{lemma}
Assume that $2<q<p<2^*$ and $(b_0)$ holds. Then, the Nehari manifold  $\mathcal{N}_\lambda$ is empty for any $\lambda \in (-\infty, \lambda_*)$.
\end{lemma}
\begin{proof}
It follows directly from item $(2)$ of Proposition \ref{projecao}.
\end{proof}

Next we prove that the sets $\mathcal{N}_\lambda^-$ and $\mathcal{N}_\lambda^+$ are $C^1$ submanifolds of $D$. This result is fundamental for the analysis of the minimization Problem \eqref{mp}.

\begin{proposition}\label{naovazia}
Suppose that $2<q<p<2^*$ and $(b_0)$ holds.  Then:
\begin{enumerate}
    \item [$(i)$] $\mathcal{N}^+_\lambda$ is not empty for every $\lambda>\lambda_*$;
    \item[$(ii)$] $\mathcal{N}^-_\lambda$ is not empty for every $\lambda_*<\lambda<\lambda_1(\Omega)$.
\end{enumerate}
Furthermore, $\mathcal{N}^+_\lambda$ and $\mathcal{N}^-_\lambda$ are $C^1$ submanifolds in $D$.
\end{proposition}
\begin{proof} 
Note that Proposition \ref{projecao} ensures that the sets $\mathcal{N}^+_\lambda$ and $\mathcal{N}^-_\lambda$ are nonempty.
To prove the second part, we define the $C^1$ functional $J_{\lambda} : D \to \mathbb{R}$ given by $J_{\lambda}(u) = E'_\lambda(u)u$, $u \in D$. It is easy to verify that $J'_\lambda(u)\phi = E''_\lambda(u)(u, \phi) + E'_\lambda(u)\phi$, $u \in D$ and $\phi\in H_0^1(\Omega)$. In particular, $J'_\lambda(u)u = E''_\lambda(u)u^2 + E'_\lambda(u)u \neq 0$ for every $u \in \mathcal{N}_\lambda^{\pm}$. Hence, the desired result follows from the Implicit Function Theorem \cite{drabek}.
\end{proof}
\begin{remark}
The non-emptiness of \(\mathcal N_\lambda^-\) is used below only for \(\lambda\in(\lambda_*,\lambda_1(\Omega))\), as established in Proposition naovazia. No global non-emptiness assertion for \(\mathcal N_\lambda^-\) is required.
\end{remark}

{
\begin{lemma}
Assume that $2<q<p<2^*$ and $(b_0)$ holds. If $u\in\mathcal{N}_\lambda^0$, then
$$||u||_{p,b}^p=\dfrac{(q-2)}{(p-q)}(||u||^2-\lambda ||u||^2_2), \ \ 
||u||_{q,b}^q=\dfrac{(p-2)}{(p-q)}(||u||^2-\lambda ||u||^2_2) \mbox{ and } ||u||_{p,b}^p=\dfrac{(q-2)}{(p-2)}||u||_{q,b}^q.$$
\end{lemma}
\begin{proof}
The identities are a direct consequence of $E'_\lambda(u)u=0$ and $E''_\lambda(u)u^2=0$.
\end{proof}
}

Next, we present a result characterizing the energy functional restricted to the Nehari manifold. 

\begin{lemma}\label{importante}
Assume that $2<q<p<2^*$ and $(b_0)$ holds.
\begin{enumerate}
\item  If $ u\in \mathcal{N}_\lambda$, then
\begin{eqnarray}
\label{functional2.0}E_{\lambda}(u)
  &=& \left(\dfrac{1}{2}-\dfrac{1}{q}\right)\left(||u||^2-\lambda||u||^2_2\right)+\left(\dfrac{1}{p}-\dfrac{1}{q}\right)||u||_{p,b}^p.
\end{eqnarray}
\item If $u\in \mathcal{N}_\lambda^0,$ then
    \begin{eqnarray}
    \label{functional2} E_{\lambda}(u) &=&\dfrac{(p-2)(q-2)}{2pq}(||u||^2-\lambda ||u||^2_2).
\end{eqnarray}
    \item If $u\in \mathcal{N}_\lambda,$ then
    \begin{eqnarray}\label{2derivada}
        E''_\lambda(u)u^2
        &=&(p-2)||u||_{p,b}^p-(q-2)||u||_{q,b}^q\nonumber\\[2ex]
    &=&(2-q)(||u||^2-\lambda||u||^2_2)+\left(p-q\right)||u||_{p,b}^p\nonumber\\[2ex]
    &=&(2-p)(||u||^2-\lambda||u||^2_2)+\left(p-q\right)||u||_{q,b}^q.
    \end{eqnarray}
    \end{enumerate} 
\end{lemma}
\begin{proof}
The proofs of the above results are based on simple calculations. Therefore, they will be omitted here.
\end{proof}

 In order to solve the minimization Problem \eqref{mp}, it is crucial to prove that the minimizer lies within the Nehari manifold. The subsequent lemma ensures this property.
\begin{lemma}\label{bel-D}
    If $\lambda_*<\lambda\leq\lambda_1(\Omega)$, then $\mathcal{N}_\lambda\subset D$.
    \end{lemma}
    \begin{proof}
     If $u\not\in D$, then $$\lambda<\lambda_1(\Omega) \leq R_n(tu)=\frac{\|u\|^2}{\|u\|_2^2},~\text{for all} ~ t>0.$$
     This implies that the equation $R_n(tu)=\lambda$ does not have {a} solution, i.e., $u \notin \mathcal{N}_\lambda$.

 If  $\lambda = \lambda_1(\Omega)$, since $\varphi_1 \in D$, for any $u \notin D$,  we have  $R_n(tu) = \frac{\|u\|^2}{\|u\|_2^2} > \lambda_1(\Omega)$. Consequently, the equation $R_n(tu) = \lambda_1(\Omega)$ has no solution.
Finally, item (2) of Proposition \ref{projecao} implies that the fibering map $t \mapsto E_{\lambda}(tu)$ has no critical point whenever $\lambda < \lambda_n(u)$.
\end{proof}

%
 It is expected that minimizers restricted to the Nehari manifold are solutions to Problem \eqref{eq1}. The following lemma ensures this.
\begin{lemma}\label{min-sol}
Suppose $2<q<p<2^*$ and {that} $(b_0)$ {holds}. Let $u \in D$ be a local minimum for $E_{\lambda}$ on $\mathcal{N}_{\lambda}^+$ or $\mathcal{N}_{\lambda}^-$, that is, assume that $u \in \mathcal{N}_{\lambda}^- \cup \mathcal{N}_{\lambda}^+$ is a minimizer for $E_{\lambda}$ restricted to $\mathcal{N}_{\lambda}^+$ or $\mathcal{N}_{\lambda}^-$, respectively. Then $u$ is a critical point of $E_{\lambda}$ on $H_0^1(\Omega)$, that is, we obtain that $E'_{\lambda}(u)\psi = 0$ for each $\psi \in H_0^1(\Omega)$.
\end{lemma}
\begin{proof}
The proof follows using a standard minimization argument in the Nehari set $\mathcal{N}_{\lambda}$. Let $u \in D$ be a local minimum for $E_{\lambda}$ in $\mathcal{N}_{\lambda}^+ \cup \mathcal{N}_{\lambda}^-$. It follows from the Lagrange Multipliers Theorem \cite{drabek} that there exists $\theta \in \mathbb{R}$ such that
$
E'_{\lambda}(u)\psi = \theta \mathcal{P}'(u)\psi, \text{for } \psi \in H_0^1(\Omega),
$
where $\mathcal{P} : D  \to \mathbb{R}$ is given by $\mathcal{P}(u) = E'_{\lambda}(u)u$. Since $0\neq E''_{\lambda}(u)u^2 = \mathcal{P}'(u)u  $, we obtain that $\theta = 0$ and $E'_{\lambda}(u)\psi = 0$ holds for each $\psi \in H_0^1(\Omega)$. This ends the proof.
\end{proof}

To establish the existence of a solution to the Problem \eqref{mp}, it is first necessary to prove that the minimizing sequence is bounded and that the corresponding minimizers are nontrivial. The following lemma provides a fundamental tool for achieving these results.
 \begin{lemma}\label{coercive}
      Let $2<q<p<2^*,$ $(b_0)$\ and\ $\ \lambda_*<\lambda<\lambda_1(\Omega)$ hold. 
      If $u  \in \mathcal{N}_\lambda$, 
    then: 
     
      \begin{enumerate}
          \item[$(1)$] {for each \(\lambda\in(\lambda_*,\lambda_1(\Omega))\), there exists a constant \(\delta_\lambda>0\) such that \(\|u\|\geq\delta_\lambda\) for every \(u\in\mathcal N_\lambda\). Moreover, the estimate is uniform with respect to \(\lambda\) on compact subintervals: for every compact interval \(K\subset\subset(\lambda_*,\lambda_1(\Omega))\), there exists a constant \(\delta_K>0\), independent of \(\lambda\in K\), such that \(\|u\|\geq\delta_K\) for every \(\lambda\in K\) and every \(u\in\mathcal N_\lambda\).
          \item[$(2)$] If \(u\in\mathcal N_\lambda^-\cup\mathcal N_\lambda^0\), then there exists a constant \(A_{\delta_\lambda}>0\) such that \(E_\lambda(u)\geq A_{\delta_\lambda}\). Moreover, the estimate is uniform with respect to \(\lambda\) on compact subintervals: for every compact interval \(K\subset\subset(\lambda_*,\lambda_1(\Omega))\), there exists a constant \(A_K>0\), independent of \(\lambda\in K\), such that \(E_\lambda(u)\geq A_K\) for every \(\lambda\in K\) and every \(u\in\mathcal N_\lambda^-\cup\mathcal N_\lambda^0\)}.
          \item[$(3)$]  \( E_\lambda \) is coercive.  
        \end{enumerate}
          
\end{lemma}
      \begin{proof}
         (1)\ Since \( u \in \mathcal{N}_\lambda \), by using \eqref{nehari} and the fact that $\lambda<\lambda_1(\Omega)$, there exists a positive constant $C=C(\lambda)=\frac{\lambda_1(\Omega)-\lambda}{\lambda_1(\Omega)}$  such that 
          \begin{eqnarray}\label{c1}
              C||u||^2\leq ||u||^2-\lambda||u||^2_2= -||u||^p_{p,b}+||u||^q_{q,b}\leq ||u||^q_{q,b}.
          \end{eqnarray}
Now, in view of \eqref{imersao1}, we obtain
\begin{eqnarray*}
     C||u||^2\leq \dfrac{||b||_\infty}{S_q^q}||u||^q.
\end{eqnarray*}
Consequently,
\begin{eqnarray}\label{delta}
    ||u||\geq \left[\dfrac{CS_q^q}{||b||_\infty}\right]^{\frac{1}{q-2}}:=\delta_\lambda>0.
    \end{eqnarray}
    
    {If $K\subset\subset(\lambda_*,\lambda_1(\Omega))$ is compact, then
    \[
    C=1-\frac{\lambda}{\lambda_1(\Omega)}
    \geq 1-\frac{\max K}{\lambda_1(\Omega)}=:C_K>0,
    \qquad\text{for every }\lambda\in K.
    \]
    Hence \eqref{delta}yields the uniform lower bound
    \[
    \|u\|\geq
    \left[\frac{C_KS_q^q}{\|b\|_\infty}\right]^{\frac1{q-2}}
    =:\delta_K>0,
    \]
    which establishes the uniform estimate in item (1).}\\
        (2) Let $u\in \mathcal{N}_\lambda^-\cup\mathcal{N}_\lambda^0. $ Then, in view of \eqref{2derivada}, we have
    $$\dfrac{q-2}{p-q}||u||^2_\lambda:=\dfrac{q-2}{p-q}(\|u\|^2-\lambda\|u\|_2^2)\geq ||u||_{p,b}^p.$$ As a consequence of \eqref{functional2.0}, the inequality above, \eqref{c1} and \eqref{delta}, we can conclude that $$E_\lambda(u)\geq  \dfrac{(p-2)(q-2)}{2pq}||u||_\lambda^2\geq \dfrac{(p-2)(q-2)}{2pq}C||u||^2\geq \dfrac{(p-2)(q-2)}{2pq}c_1\delta^2_\lambda:=A_{\delta_\lambda}.$$  {If $\lambda\in K\subset\subset(\lambda_*,\lambda_1(\Omega))$, the preceding estimates yield
    \[
    E_\lambda(u)\geq
    \frac{(p-2)(q-2)}{2pq}\,C_K\delta_K^2
    =:A_K>0,
    \]
    for every \(\lambda\in K\) and every \(u\in\mathcal N_\lambda^-\cup\mathcal N_\lambda^0\). Since \(A_K\) is independent of \(\lambda\in K\), this establishes the uniform estimate in item (2).}\\
  (3)
  Since
   $$E_\lambda(u)=\dfrac{1}{2}||u||_\lambda^2+\dfrac{1}{p}||u||^p_{p,b}-\dfrac{1}{q}||u||^q_{q,b}$$ 
  and $q<p$, by \eqref{desigualdade} 
  we obtain the existence of a positive constant  $C=(|\Omega|\|b\|_\infty)^{\frac{p-q}{p}}$ such that
 \begin{eqnarray}\label{coe-lim}
     E_\lambda(u)\geq \dfrac{1}{2}||u||_\lambda^2+\dfrac{1}{p}||u||^p_{p,b}-\dfrac{C}{q}||u||^q_{p,b}=\dfrac{1}{2}||u||_\lambda^2+||u||^p_{p,b}\left(\dfrac{1}{p}-\dfrac{C}{q}||u||^{q-p}_{p,b}\right).
 \end{eqnarray}
  Now, suppose that $||u||_\lambda\to \infty$ and that there exists $C>0$ such that $||u||_{p,b}\leq C$. Then $E_\lambda(u)\to \infty.$ On the other hand, if $||u||_{p,b}\to \infty,$ we again obtain $E_\lambda(u)\to \infty.$
    \end{proof}
{
\begin{lemma}\label{bounded-comp}  Assume $2<q<p<2^*,$ $(b_0)$\ and\ $\lambda_1(\Omega)\le\lambda<{\lambda_1(\Omega_0)}$ hold. 
     Let $(u_k)\subset \mathcal{N_\lambda^{+}}$ be a minimizing sequence. Then $(u_k)$ is bounded.
\begin{proof}
{
Suppose, by contradiction, that $\|u_k\|\to\infty$, and set
$v_k=u_k/\|u_k\|$. Then $\|v_k\|=1$ and, up to a subsequence,
$v_k\rightharpoonup v$ in $H_0^1(\Omega)$ and
$v_k\to v$ in $L^r(\Omega)$ for every $2\le r<2^*$.

Since $u_k\in\mathcal N_\lambda$, the Nehari identity yields
\[
\|u_k\|^2-\lambda\|u_k\|_2^2
 =\int_\Omega b(x)\bigl(|u_k|^q-|u_k|^p\bigr)\,dx.
\]
Since the function $t\mapsto t^q-t^p$ is bounded above on $[0,\infty)$; there exists a constant $C>0$, independent of $k$, such that
\(
\|u_k\|^2-\lambda\|u_k\|_2^2\le C.
\)
Dividing by \(\|u_k\|^2\), we obtain
$
1-\lambda\|v_k\|_2^2
\le {C}/{\|u_k\|^2}.
$
Since \(\|u_k\|\to\infty\) and \(v_k\to v\) strongly in \(L^2(\Omega)\), passing to the limit yields
$
1\le\lambda\|v\|_2^2,
$
which implies that \(v\neq0\).

On the other hand, dividing the Nehari identity by $\|u_k\|^p$ we obtain
\[
\|u_k\|^{2-p}\bigl(1-\lambda\|v_k\|_2^2\bigr)
=\|u_k\|^{q-p}\|v_k\|_{q,b}^q-\|v_k\|_{p,b}^p.
\]
Since $p>q>2$, $\|u_k\|\to\infty$, and $(v_k)$ is bounded in
$H_0^1(\Omega)$, both
$\|u_k\|^{2-p}(1-\lambda\|v_k\|_2^2)$ and
$\|u_k\|^{q-p}\|v_k\|_{q,b}^q$ converge to zero. Consequently,
$
\|v_k\|_{p,b}^p\longrightarrow0.
$
Since \(v_k\to v\) strongly in \(L^p(\Omega)\), it follows that \(\|v\|_{p,b}=0\). By assumption \((b_0)\), we have \(v=0\) almost everywhere in \(\Omega\setminus\Omega_0\), and hence \(v\in H_0^1(\Omega_0)\).

Finally, by weak lower semicontinuity, the variational characterization of \(\lambda_1(\Omega_0)\), and the previously established inequality \(1\leq\lambda\|v\|_2^2\), we obtain
\[
1\ge \|v\|^2
\ge \lambda_1(\Omega_0)\|v\|_2^2
\ge \frac{\lambda_1(\Omega_0)}{\lambda}>1,
\]
which is a contradiction. Thus $(u_k)$ is bounded.
}
\end{proof}
\end{lemma}
}
Finding a minimizer of Problem~\eqref{mp} requires an understanding of the closures of the sets $\mathcal{N}_\lambda^{\pm}$.
\begin{lemma}\label{fecho} Suppose $2<q<p<2^*$, $(b_0)$ and $\lambda_*<\lambda<\lambda_1(\Omega)$. Then,
\(
\overline{\mathcal{N}_{\lambda}^{\pm}}
\subset
\mathcal{N}_{\lambda}^{\pm}\cup\mathcal{N}_{\lambda}^{0}.
\)
\end{lemma}
\begin{proof}
We prove the assertion for $\mathcal{N}_\lambda^+$; the proof for
$\mathcal{N}_\lambda^-$ is analogous. Let $(u_k)\subset\mathcal{N}_\lambda^+$
and suppose that $u_k\to u$ in $H_0^1(\Omega)$. Then
\(
E_\lambda'(u_k)u_k=0
\  \text{and}\  
E_\lambda''(u_k)u_k^2>0.
\)
By the continuity of $E_\lambda'$ and $E_\lambda''$, we get
\(
E_\lambda'(u)u=0
\  \text{and}\ 
E_\lambda''(u)u^2\geq0.
\)
Moreover, item $(1)$ of Lemma~\ref{coercive} gives
$\|u_k\|\geq\delta_\lambda>0$, and hence $u\neq0$. Therefore
\(
u\in\mathcal{N}_\lambda^+\cup\mathcal{N}_\lambda^0
\), which proves the desired inclusion.
\end{proof}

\section{Existence and Multiplicity of Solutions} \label{SB3}

{The following result shows that, under assumption $(b_0)$ alone, every nonnegative nontrivial solution is in fact positive.}

\begin{proposition}\label{pos-sol}
{Suppose that $2<q<p<2^*$ and assumption $(b_0)$ holds}. Let $w$ be a nonnegative, nontrivial weak solution of Problem~\eqref{eq1}. Then $w>0$ in $\Omega$.
\end{proposition}
\begin{proof}
{By standard elliptic regularity and bootstrap arguments} for subcritical semilinear equations (see, e.g., \cite{drabek}), $w\in L^\infty(\Omega)$. Since $\Omega$ is smooth and $w=0$ on $\partial\Omega$, {elliptic regularity further yields} $w\in W^{2,r}(\Omega)$, for every finite $r$. In particular, $w\in C^{1,\beta}(\overline\Omega)$ for some $\beta\in(0,1)$. Consequently
\[
a(x):=\lambda-b(x)\big(w^{p-2}-w^{q-2}\big)\in L^\infty(\Omega),
\]
and $w$ satisfies
$
-\Delta w=a(x)w~\text{in }\Omega.
$
The strong maximum principle for {nonnegative} weak solutions with bounded zero-order coefficient (see, e.g., \cite{pucci}) implies that either $w\equiv0$ or $w>0$ in $\Omega$. Since $w$ is nontrivial, {we conclude that} $w>0$ in $\Omega$.
\end{proof}

\subsection{Solution in $\mathcal{N}_\lambda^+$}

{
This subsection is devoted to establishing the existence of a minimizer for the  problem
\[
c_\lambda^+ := \inf \{ E_{\lambda}(u) : u \in \mathcal{N}^{+}_{\lambda} \}.
\]
{Beyond establishing the existence of such a minimizer throughout the parameter interval considered below, this result plays a crucial role in proving the existence of a second solution, which will be addressed in the following section.}
}
{
\begin{lemma}\label{solutionN+}
    Suppose that $2<q<p<2^*$, assumption $(b_0)$ holds, and
$\lambda_*<\lambda<\lambda_1(\Omega_0)$.
If $(u_k)\subset\mathcal{N}_\lambda^+$ is a minimizing sequence
for $E_\lambda$ on $\mathcal{N}_\lambda^+$, then there exists $
u\in(\mathcal{N}_\lambda^+\cup\mathcal{N}_\lambda^0)
\setminus\{0\}$
such that $u_k\to u$ strongly in $H_0^1(\Omega)$.
\end{lemma}
\begin{proof}
If $\lambda<\lambda_1(\Omega)$, the sequence $(u_k)$ is bounded
in $H_0^1(\Omega)$ by Lemma~\ref{coercive}. If
$\lambda\ge\lambda_1(\Omega)$, its boundedness follows from
Lemma~\ref{bounded-comp}. Thus, $(u_k)$ is bounded in
$H_0^1(\Omega)$. Passing to a subsequence, we obtain
$
u_k\rightharpoonup u~\text{in }H_0^1(\Omega),
~
u_k\to u~\text{in }L^r(\Omega),
~ 2\le r<2^*.
$
In particular, since $b\in L^\infty(\Omega)$, we have
\[
\int_\Omega b(x)|u_k|^p\,dx
\to \int_\Omega b(x)|u|^p\,dx
\ \text{ and }\ 
\int_\Omega b(x)|u_k|^q\,dx
\to \int_\Omega b(x)|u|^q\,dx.
\]

We first show that $u\ne0$. If $\lambda<\lambda_1(\Omega)$ and
$u=0$, then the Nehari identity, together with the compact
embeddings, implies that $\|u_k\|\to0$, contradicting
Lemma~\ref{coercive}, item~$(1)$.

If $\lambda\ge\lambda_1(\Omega)$, Corollary~\ref{energ-signal},
applied to $\varphi_1$, yields an element of $\mathcal N_\lambda^+$
with negative energy. Consequently, $c_\lambda^+<0$. If $u=0$,
the Nehari identity again implies that $\|u_k\|\to0$, and hence
$E_\lambda(u_k)\to0$. This contradicts
$E_\lambda(u_k)\to c_\lambda^+<0$. Therefore, $u\ne0$.

We now prove that $u_k\to u$ strongly in $H_0^1(\Omega)$. Suppose,
by contradiction, that this convergence does not hold. Passing
to a subsequence if necessary, there exists $\varepsilon_0>0$
such that $
\|u_k-u\|\ge\varepsilon_0
~\text{for every }k.
$
Since $H_0^1(\Omega)$ is a Hilbert space and
$u_k\rightharpoonup u$, we have $
\|u_k\|^2
=
\|u\|^2+\|u_k-u\|^2+o_k(1).
$
In particular,
$
\liminf_{k\to\infty}\|u_k\|^2
\ge \|u\|^2+\varepsilon_0^2.
$
Passing to a further subsequence, we may assume that $
\|u_k-u\|^2\to d
$
for some $d>0$. By the compact embeddings,
$
u_k\to u\quad\text{in }L^r(\Omega),
~2\le r<2^*.
$
Consequently, the $L^2$-term and the weighted $p$- and $q$-terms
in $E_\lambda(u_k)$ converge. It follows that
\[
E_\lambda(u_k)
=
E_\lambda(u)+\frac12\|u_k-u\|^2+o_k(1),
\]
and hence
$
E_\lambda(u)
=
c_\lambda^+- d/2
<c_\lambda^+.
$

Similarly, since $E_\lambda'(u_k)u_k=0$, we obtain
\[
0
=
E_\lambda'(u_k)u_k
=
E_\lambda'(u)u+\|u_k-u\|^2+o_k(1),
\]
which yields
$
E_\lambda'(u)u=-d<0.
$

If $u\notin D$, then $u\in H_0^1(\Omega_0)$, and therefore
\[
\begin{aligned}
E_\lambda'(u)u
&=\|u\|^2-\lambda\|u\|_2^2\ge
\bigl(\lambda_1(\Omega_0)-\lambda\bigr)\|u\|_2^2
>0,
\end{aligned}
\]
a contradiction. Thus $u\in D$.

Since $E_\lambda'(u)u<0$, we have $R_n(u)<\lambda$.
Consequently, the largest Nehari projection
$t_{\lambda,+}(u)$ exists and satisfies $t_{\lambda,+}(u)>1$.
Moreover, $R_n(tu)<\lambda$ for
$1\le t<t_{\lambda,+}(u)$, so the map
$t\mapsto E_\lambda(tu)$ is strictly decreasing on this interval.
By continuity, we conclude that
$
E_\lambda(t_{\lambda,+}(u)u)
<
E_\lambda(u)
<
c_\lambda^+.
$
This contradicts the definition of $c_\lambda^+$, since
$t_{\lambda,+}(u)u\in\mathcal N_\lambda^+$.
Therefore,
$
u_k\to u\quad\text{strongly in }H_0^1(\Omega).
$

Finally, since $E_\lambda'(u_k)u_k=0$ and
$E_\lambda''(u_k)u_k^2>0$, the strong convergence just established,
together with the continuity of $E_\lambda'$ and $E_\lambda''$,
implies that
$
E_\lambda'(u)u=0,
\qquad
E_\lambda''(u)u^2\ge0.
$
Hence $u\in\mathcal N_\lambda^+\cup\mathcal N_\lambda^0$.
Moreover, $u\ne0$ by the preceding argument.
\end{proof}
}

For the minimizer to be a solution, it must belong to $D$. The following proposition establishes this fact.
\begin{proposition}\label{belo-D}
    Suppose $2<q<p<2^*$, $(b_0)$ and  $\lambda>\lambda_*$. Let $u\in \mathcal{N}_\lambda^+\cup \mathcal{N}_\lambda^0$ be a function such that $c_\lambda^+=E_\lambda(u)$. Then $u\in D$.
\end{proposition}
\begin{proof}
    If $\lambda_*<\lambda\leq  \lambda_1(\Omega)$, by Lemma \ref{bel-D} we have {$u\in D$}.

    We now consider the case $\lambda >\lambda_1(\Omega)$. Assume, by contradiction, that $u\notin D$. Then $\lambda=R_n(u)=\frac{\|u\|^2}{\|u\|_2^2}$, which implies $E_\lambda(u)=0$. On the other hand,
$E_\lambda(u)=c_\lambda^+\leq E_\lambda(t_{\lambda,+}(\varphi_1)\varphi_1)<0,$
   which yields a contradiction.    
\end{proof}

{ Our goal is to show that a minimizer on the Nehari manifold is a solution to Problem \eqref{eq1}. By Lemma \ref{min-sol}, it suffices to prove that this minimizer does not belong to $\mathcal{N}_\lambda^0$. The main difficulty arises when
$\lambda\in(\bar{\lambda},\lambda^*),$ where $\bar{\lambda}$ was defined in \eqref{lambdane}. The following lemma will be essential in addressing this case.} 
{
\begin{lemma}\label{c+cont}
One has $
c_{\lambda^*}^+=0
~\text{and}~
\lim_{\lambda\uparrow\lambda^*}c_\lambda^+=0.
$
\end{lemma}

\begin{proof}
Let $u_{\lambda^*}=s(u^*)u^*$ be the zero-energy solution
given by Proposition~\ref{le}. Since
$s(u^*)>\tau(u^*)$ and
\(
R_n(u_{\lambda^*})=R_e(u_{\lambda^*})=\lambda^*,
\)
Proposition~\ref{1} implies that
$u_{\lambda^*}\in\mathcal N_{\lambda^*}^+$.

By the definition of $\lambda^*$,
$
R_e(v)\ge\lambda^*
~\text{for every }v\in D.
$
Consequently, $E_{\lambda^*}(v)\ge0$ for every $v\in D$.
If $v\notin D$, then
$
E_{\lambda^*}(v)
=1/2\bigl(\|v\|^2-\lambda^*\|v\|_2^2\bigr)>0,
$
since $\lambda^*<\lambda_1(\Omega)$. Thus
$E_{\lambda^*}(v)\ge0$ for every $v\in H_0^1(\Omega)$.
Since $u_{\lambda^*}\in\mathcal N_{\lambda^*}^+$ and
$E_{\lambda^*}(u_{\lambda^*})=0$, we conclude that
$c_{\lambda^*}^+=0$.

Now let $\lambda<\lambda^*$. The same argument shows that
$E_\lambda(v)>0$ for every $v\ne0$, and hence
$c_\lambda^+\ge0$. Moreover,
$
\lambda_n(u_{\lambda^*})
<
R_n(u_{\lambda^*})
=\lambda^*.
$
Therefore, for $\lambda<\lambda^*$ sufficiently close to
$\lambda^*$, the projection
$t_{\lambda,+}(u_{\lambda^*})$ is well defined. By
Proposition~\ref{projecao},
$
t_{\lambda,+}(u_{\lambda^*})\longrightarrow1
~\text{as }\lambda\uparrow\lambda^*.
$
It follows that
\[
\begin{aligned}
0\le c_\lambda^+
&\le
E_\lambda\bigl(t_{\lambda,+}(u_{\lambda^*})
u_{\lambda^*}\bigr)\longrightarrow E_{\lambda^*}(u_{\lambda^*})=0
\qquad\text{as }\lambda\uparrow\lambda^*.
\end{aligned}
\]
The desired conclusion follows.
\end{proof}

{ Next, we will find a range of parameters for which the minimizers on the Nehari manifold belong to $\mathcal{N}_\lambda^+$ which will guarantee us that this minimizer is a solution to the Problem \eqref{eq1}.}
\begin{proposition}\label{sol-N+}
Assume that $2<q<p<2^*$ and that $(b_0)$ holds. Then there
exists $\zeta>0$ such that, for every $
\lambda\in(\lambda^*-\zeta,\lambda_1(\Omega_0)),
$
there exists $u_\lambda\in\mathcal N_\lambda^+$ satisfying
$
c_\lambda^+=E_\lambda(u_\lambda).
$
\end{proposition}

\begin{proof}
We consider three cases.

\smallskip
\noindent
\textit{Case 1: $\lambda^*-\zeta<\lambda<\lambda^*$.}
Choose $\zeta_0>0$ sufficiently small that $
K=[\lambda^*-\zeta_0,\lambda^*]
\subset\subset(\lambda_*,\lambda_1(\Omega)).
$
By item~$(2)$ of Lemma~\ref{coercive}, there exists $A_K>0$,
independent of $\lambda\in K$, such that
$
E_\lambda(u)\ge A_K
~\text{for every }
u\in\mathcal N_\lambda^-\cup\mathcal N_\lambda^0.
$
By Lemmas~\ref{fecho} and~\ref{solutionN+}, the infimum
$c_\lambda^+$ is attained at some
$
u_\lambda\in\mathcal N_\lambda^+\cup\mathcal N_\lambda^0.
$
Moreover, Lemma~\ref{c+cont} gives
$c_\lambda^+\to0$ as $\lambda\uparrow\lambda^*$.
Shrinking $\zeta_0$ if necessary, we may therefore assume that
$
0\le c_\lambda^+<A_K
~\text{for }\lambda^*-\zeta_0<\lambda<\lambda^*.
$
It follows that $u_\lambda\notin\mathcal N_\lambda^0$, and hence
$u_\lambda\in\mathcal N_\lambda^+$.

\smallskip
\noindent
\textit{Case 2: $\lambda=\lambda^*$.}
By Proposition~\ref{le}, there exists a zero-energy solution
$u_{\lambda^*}\in\mathcal N_{\lambda^*}^+$. Since
$c_{\lambda^*}^+=0$, we have
$
E_{\lambda^*}(u_{\lambda^*})=c_{\lambda^*}^+.
$

\smallskip
\noindent
\textit{Case 3: $\lambda^*<\lambda<\lambda_1(\Omega_0)$.}
By Lemma~\ref{solutionN+}, the infimum $c_\lambda^+$ is attained
at some
$
u\in\mathcal N_\lambda^+\cup\mathcal N_\lambda^0.
$
Proposition~\ref{belo-D} yields $u\in D$. Since
$\lambda_e(u_{\lambda^*})=\lambda^*$, Corollary~\ref{energ-signal}
implies that
$
c_\lambda^+
\le E_\lambda\bigl(t_{\lambda,+}(u_{\lambda^*})u_{\lambda^*}\bigr)
<0.
$
On the other hand, if $u\in\mathcal N_\lambda^0\cap D$, then
\[
\|u\|_{p,b}^p
=\frac{q-2}{p-q}
\bigl(\|u\|^2-\lambda\|u\|_2^2\bigr)>0.
\]
Consequently, by \eqref{functional2},
\[
E_\lambda(u)
=\frac{(p-2)(q-2)}{2pq}
\bigl(\|u\|^2-\lambda\|u\|_2^2\bigr)>0,
\]
contradicting $E_\lambda(u)=c_\lambda^+<0$.
Thus $u\in\mathcal N_\lambda^+$.
Taking $\zeta=\zeta_0$ completes the proof.
\end{proof}

{
We next investigate the dependence of the minimizing level $c_\lambda^+$ and of the corresponding minimizers on the parameter $\lambda$.

\begin{proposition}\label{cont-cplus}
Assume $2<q<p<2^*$ and $(b_0)$. Then the map
\(
\lambda\mapsto c_\lambda^+
\)
is continuous and strictly decreasing on
\(
(\lambda^*-\zeta,\lambda_1(\Omega_0)).
\)
Moreover, if $\lambda_n\to\lambda$ and
\(
u_n\in\mathcal N_{\lambda_n}^+,
~
E_{\lambda_n}(u_n)=c_{\lambda_n}^+,
\)
then, up to a subsequence,
\(
u_n\to u
~\mbox{in }H_0^1(\Omega),
\)
where
\(
u\in\mathcal N_\lambda^+
\)
and
\(
E_\lambda(u)=c_\lambda^+.
\)
\end{proposition}
\begin{proof}
Let $\mu<\lambda$ and let $u_\mu\in\mathcal N_\mu^+$ be a minimizer. Since
\(
R_n(u_\mu)=\mu
\)
and
\(
\tau(u_\mu)<1,
\)
there exists
\(
t_{\lambda,+}(u_\mu)>1
\)
such that
\(
t_{\lambda,+}(u_\mu)u_\mu\in\mathcal N_\lambda^+.
\)
Moreover,
\(
R_n(tu_\mu)<\lambda
\)
for
\(
1\le t<t_{\lambda,+}(u_\mu).
\)
Hence
\[
c_\lambda^+
\le
E_\lambda(t_{\lambda,+}(u_\mu)u_\mu)
<
E_\lambda(u_\mu)
<
E_\mu(u_\mu)
=
c_\mu^+,
\]
and therefore $c_\lambda^+$ is strictly decreasing.

Fix now $\lambda$ and let $u_\lambda\in\mathcal N_\lambda^+$ be a minimizer. Since
\[
\frac{d}{dt}R_n(tu_\lambda)\Big|_{t=1}>0,
\]
the implicit function theorem gives
\(
t_{\mu,+}(u_\lambda)\to1
\)
as
\(
\mu\to\lambda.
\)
Thus
\[
\limsup_{\mu\to\lambda}c_\mu^+
\le
\lim_{\mu\to\lambda}
E_\mu(t_{\mu,+}(u_\lambda)u_\lambda)
=
c_\lambda^+.
\]

Let $\lambda_n\to\lambda$ and let
\(
u_n\in\mathcal N_{\lambda_n}^+
\)
be minimizers. We claim that $(u_n)$ is bounded. Otherwise, after passing to a subsequence,
\(
\|u_n\|\to\infty.
\)
Set
\(
w_n=u_n/\|u_n\|.
\)
Then, up to a subsequence,
\(
w_n\rightharpoonup w
\quad\mbox{in }H_0^1(\Omega),
~
w_n\to w
~\mbox{in }L^r(\Omega),
~ 2\le r<2^*.
\)
The Nehari identity yields
\(
1-\lambda_n\|w_n\|_2^2\le o_n(1),
\)
and hence
\(
1\le\lambda\|w\|_2^2.
\)
On the other hand, dividing the same identity by $\|u_n\|^p$ gives
\(
\|w_n\|_{p,b}^p\to0.
\)
Thus $w\in H_0^1(\Omega_0)$, and consequently
\[
1
\le
\lambda\|w\|_2^2
\le
\frac{\lambda}{\lambda_1(\Omega_0)}\|w\|^2
\le
\frac{\lambda}{\lambda_1(\Omega_0)}
<1,
\]
a contradiction.

Therefore, $(u_n)$ is bounded in $H_0^1(\Omega)$. Passing to a
subsequence, we may assume that
$
u_n\rightharpoonup u~\text{in }H_0^1(\Omega),
~
u_n\to u~\text{in }L^r(\Omega),
~ 2\le r<2^*.
$
By Proposition~\ref{belo-D} and Lemma~\ref{min-sol}, each $u_n$
is a weak solution of $(P_{\lambda_n})$. Testing its weak
formulation with $u_n-u$ and using the strong $L^r$-convergence,
we obtain
$
\|u_n-u\|\to0.
$
Consequently, by the continuity of $E_\lambda'$,
$
E_\lambda'(u)=0.
$
Moreover, $u\neq0$. Indeed, if $\lambda<\lambda_1(\Omega)$, this follows from the uniform estimate in Lemma~\ref{coercive}; if $\lambda\ge\lambda_1(\Omega)$, then $\lambda>\lambda^*$ and $c_\lambda^+<0$, while the upper semicontinuity proved above excludes $u=0$.

Since
\(
E_{\lambda_n}''(u_n)u_n^2>0,
\)
we obtain
\(
u\in\mathcal N_\lambda^+\cup\mathcal N_\lambda^0.
\)
We claim that
\(
c_\lambda^+
<
\inf_{v\in\mathcal N_\lambda^0}E_\lambda(v)
\)
for every
\(
\lambda\in(\lambda^*-\zeta,\lambda_1(\Omega_0)).
\)

Indeed, if
\(
\lambda^*-\zeta<\lambda<\lambda^*,
\)
this is precisely the separation used in Proposition~\ref{sol-N+}, at
\(
\lambda=\lambda^*,
\)
we have
\(
c_{\lambda^*}^+=0,
\)
while
\(
E_{\lambda^*}(v)>0
\)
for every
\(
v\in\mathcal N_{\lambda^*}^0.
\)
Finally, if
\(
\lambda^*<\lambda<\lambda_1(\Omega_0),
\)
then
\(
c_\lambda^+<0.
\)
Moreover, if
\(
v\in\mathcal N_\lambda^0,
\)
then necessarily
\(
v\in D.
\)
Indeed, if
\(
v\notin D,
\)
then
\(
v\in H_0^1(\Omega_0)
\)
and the Nehari identity would give
$
\|v\|^2=\lambda\|v\|_2^2,
$
which contradicts
\(
\lambda<\lambda_1(\Omega_0).
\)
Hence, by the Nehari identities,
\(
E_\lambda(v)>0.
\)
Thus the claim follows.

Consequently,
\(
u\notin\mathcal N_\lambda^0,
\)
and hence
\(
u\in\mathcal N_\lambda^+.
\)
Finally,
\[
c_\lambda^+
\le
E_\lambda(u)
=
\liminf_{n\to\infty}c_{\lambda_n}^+
\le
\limsup_{n\to\infty}c_{\lambda_n}^+
\le
c_\lambda^+.
\]
Therefore
\(
c_{\lambda_n}^+\to c_\lambda^+,
~
E_\lambda(u)=c_\lambda^+,
\)
which proves the result.
\end{proof}
}
{
\begin{cor}\label{zero-cplus}
Assume that $2<q<p<2^*$ and that $(b_0)$ holds. Then, for
$
\lambda\in(\lambda^*-\zeta,\lambda_1(\Omega_0)),
$
we have
$
c_\lambda^+=0
~\text{if and only if}~
\lambda=\lambda^*.
$
\end{cor}

\begin{proof}
Let $\lambda^*-\zeta<\lambda<\lambda^*$, and let
$u_\lambda\in\mathcal N_\lambda^+$ be a minimizer of
$E_\lambda$ on $\mathcal N_\lambda^+$. By Proposition~\ref{sol-N+},
$c_\lambda^+=E_\lambda(u_\lambda)$. Since $u_\lambda\in D$ and
$R_e(u_\lambda)\ge\lambda^*>\lambda$, we obtain
\[
c_\lambda^+
=E_\lambda(u_\lambda)
=\frac12\|u_\lambda\|_2^2
\bigl(R_e(u_\lambda)-\lambda\bigr)>0.
\]
Moreover, Lemma~\ref{c+cont} gives $c_{\lambda^*}^+=0$, whereas
Proposition~\ref{sol-N+} yields
$
c_\lambda^+<0
~\text{for }
\lambda^*<\lambda<\lambda_1(\Omega_0).
$
Therefore, $c_\lambda^+=0$ if and only if $\lambda=\lambda^*$.
\end{proof}
}

\subsection{A Mountain Pass type solution}
In this subsection, we establish the existence of a second solution to
Problem~\eqref{eq1}. Starting from the minimizing solution in
$\mathcal N_\lambda^+$ obtained above, we verify the mountain-pass
geometry of $E_\lambda$ and apply the Mountain Pass Theorem to obtain
a second critical point.
\begin{lemma}\label{geometry}
Assume that $\lambda^*-\zeta<\lambda<\lambda_1(\Omega)$.
Then the functional $E_\lambda$ satisfies the mountain-pass
geometry. More precisely, there exist
$\rho_\lambda,\eta_\lambda>0$ and
$e_\lambda\in H_0^1(\Omega)$ such that
\[
E_\lambda(u)\ge\eta_\lambda
\quad\text{if }\|u\|=\rho_\lambda,\qquad
\|e_\lambda\|>\rho_\lambda,\qquad
E_\lambda(e_\lambda)<\eta_\lambda.
\]
\end{lemma}

\begin{proof}
Since $\lambda<\lambda_1(\Omega)$, we have
\[
E_\lambda(u)\ge
\frac12\left(1-\frac{\lambda}{\lambda_1(\Omega)}\right)
\|u\|^2
-\frac{\|b\|_\infty}{qS_q^q}\|u\|^q.
\]
Thus, $E_\lambda$ is positive on a sufficiently small
sphere centered at the origin.

Suppose first that $\lambda^*\le\lambda<\lambda_1(\Omega)$.
Replacing $u_{\lambda^*}$ by $|u_{\lambda^*}|$ if necessary,
we may assume that $u_{\lambda^*}\ge0$, since this
replacement preserves the energy and the fibering map.
Set
$
e_\lambda=t_{\lambda,+}(u_{\lambda^*})u_{\lambda^*}.
$
By Corollary~\ref{energ-signal},
$
E_\lambda(e_\lambda)\le0.
$
Since $e_\lambda\neq0$, we may choose $\rho_\lambda>0$
sufficiently small so that
$\rho_\lambda<\|e_\lambda\|$ and
\[
\eta_\lambda:=
\inf_{\|u\|=\rho_\lambda}E_\lambda(u)>0.
\]
This gives the desired inequalities.

Now suppose that $\lambda^*-\zeta<\lambda<\lambda^*$.
As above, we may assume that the minimizing solution
$u_\lambda\in\mathcal N_\lambda^+$ is nonnegative,
since replacing it by its absolute value preserves both
the energy and the Nehari constraint.

Choose a compact interval
$K\subset\subset(\lambda_*,\lambda_1(\Omega))$ containing
$[\lambda^*-\zeta,\lambda^*]$. By Lemma~\ref{coercive},
the minimizing solutions satisfy
$
\|u_\lambda\|\ge\delta_K>0.
$
Moreover, choosing $\rho\in(0,\delta_K)$ sufficiently
small, the preceding estimate yields
\[
\eta:=
\inf_{\substack{\lambda\in K\\ \|u\|=\rho}}
E_\lambda(u)>0.
\]
By Lemma~\ref{c+cont}, after decreasing $\zeta$ if
necessary, we have
$
c_\lambda^+=E_\lambda(u_\lambda)<\eta ~\text{for }\lambda\in(\lambda^*-\zeta,\lambda^*).
$
Consequently, taking
$
e_\lambda=u_\lambda,~
\rho_\lambda=\rho,~
\eta_\lambda=\eta,
$
we obtain the desired mountain-pass geometry.
\end{proof}

{
For the mountain-pass construction, we fix, for each $\lambda$, a
nonnegative endpoint $e_\lambda$ satisfying the conditions of
Lemma~\ref{geometry}. In the case $\lambda>\lambda_e(\varphi_1)$,
we specifically choose
$
e_\lambda=t_{\lambda,+}(\varphi_1)\varphi_1.
$
By the definition of $\lambda_e(\varphi_1)$, this choice satisfies
$E_\lambda(e_\lambda)<0$ and $\|e_\lambda\|>\rho_\lambda$.
This particular choice will be used in the asymptotic estimate
as $\lambda\to\lambda_1(\Omega)^-$.
}

\begin{proposition}\label{ps}
Assume that $\lambda^*-\zeta<\lambda<\lambda_1(\Omega)$.
Then the functional $E_\lambda$ satisfies the Palais--Smale
(PS) condition.
\end{proposition}

\begin{proof}
Let $(w_k)$ be a Palais--Smale sequence, that is,
\[
|E_\lambda(w_k)|\le C
\qquad\text{and}\qquad
\|E'_\lambda(w_k)\|_{H^{-1}}=o_k(1).
\]
Since $\lambda<\lambda_1(\Omega)$, estimate \eqref{coe-lim}
implies that $(w_k)$ is bounded in $H_0^1(\Omega)$.
Hence, up to a subsequence,
$
w_k\rightharpoonup w~\text{in }H_0^1(\Omega),
~
w_k\to w~\text{in }L^r(\Omega),
~ 2\le r<2^*.
$
In particular,
$
w_k\to w\quad\text{in }
L^2(\Omega)\cap L^p(\Omega)\cap L^q(\Omega).
$
Since $b\in L^\infty(\Omega)$ and $(w_k)$ is bounded in
$L^p(\Omega)\cap L^q(\Omega)$, we obtain
\[
\begin{aligned}
&\int_\Omega b(x)
\big(|w_k|^{p-2}w_k-|w_k|^{q-2}w_k\big)(w_k-w)\,dx=o_k(1).
\end{aligned}
\]
Moreover, the strong convergence in $L^2(\Omega)$ gives
\[
\lambda\int_\Omega w_k(w_k-w)\,dx=o_k(1).
\]
Testing $E'_\lambda(w_k)$ with $w_k-w$, we therefore get
\[
\int_\Omega\nabla w_k\cdot\nabla(w_k-w)\,dx=o_k(1).
\]
On the other hand, the weak convergence in $H_0^1(\Omega)$
implies that
\[
\int_\Omega\nabla w\cdot\nabla(w_k-w)\,dx=o_k(1).
\]
Subtracting these two identities, we conclude that
$
\|w_k-w\|^2=o_k(1).
$
Thus, $w_k\to w$ strongly in $H_0^1(\Omega)$, and the
Palais--Smale condition follows.
\end{proof}
{
For the fixed endpoint $e_\lambda$ chosen above, we define
$
\Gamma_\lambda
:=
\left\{
\gamma\in C([0,1],H_0^1(\Omega)):
\gamma(0)=0,\ \gamma(1)=e_\lambda
\right\}.
$
}
{
\begin{lemma}\label{mp-solution}
Let \(\lambda\in(\lambda^*-\zeta,\lambda_1(\Omega))\) and let \(e_\lambda\ge0\) be fixed as above. Then Problem~\eqref{eq1} admits a positive mountain-pass solution \(v_\lambda\) satisfying
\[
E_\lambda(v_\lambda)
=
c_\lambda^m
:=
\inf_{\gamma\in\Gamma_\lambda}
\max_{r\in[0,1]}E_\lambda(\gamma(r))
\ge\eta_\lambda>0.
\]
In particular, \(v_\lambda\in D\).
\end{lemma}

\begin{proof}
For \(u\in H_0^1(\Omega)\), set \(u^+=\max\{u,0\}\) and define
\[
\widehat E_\lambda(u)
=
\frac12\|u\|^2
-\frac{\lambda}{2}\|u^+\|_2^2
+\frac1p\|u^+\|_{p,b}^p
-\frac1q\|u^+\|_{q,b}^q.
\]
Since \(\lambda<\lambda_1(\Omega)\),
\[
\widehat E_\lambda(u)
\ge
\frac12\left(1-\frac{\lambda}{\lambda_1(\Omega)}\right)\|u\|^2
-\frac{\|b\|_\infty}{qS_q^q}\|u\|^q,
\]
and hence \(\widehat E_\lambda\) is positive on a sufficiently small sphere. Moreover, \(e_\lambda\ge0\), and therefore
\(
\widehat E_\lambda(e_\lambda)
=
E_\lambda(e_\lambda)
<
\eta_\lambda.
\)
Thus \(\widehat E_\lambda\) has the mountain-pass geometry.

Furthermore, \(\widehat E_\lambda\) satisfies the Palais--Smale condition. Indeed, any Palais--Smale sequence is bounded, since \(\lambda<\lambda_1(\Omega)\), and the compact embeddings
\(
H_0^1(\Omega)\hookrightarrow L^r(\Omega),
~ 2\le r<2^*,
\)
together with the continuity of \(u\mapsto u^+\), yield strong convergence exactly as in Proposition~\ref{ps}.

Let \(\widehat c_\lambda^m\) be the corresponding mountain-pass level. Since \(e_\lambda\ge0\), for every \(\gamma\in\Gamma_\lambda\) the path \(|\gamma|\) also belongs to \(\Gamma_\lambda\), and
\(
E_\lambda(|\gamma(r)|)=E_\lambda(\gamma(r)).
\)
Hence \(c_\lambda^m\) may be computed by restricting the admissible paths to the positive cone. On the other hand,
$
\widehat E_\lambda(u)
=
\widehat E_\lambda(u^+)+1/2\|u^-\|^2,
$
so every admissible path for \(\widehat E_\lambda\) can be replaced by its positive part without increasing its maximum energy. Since \(\widehat E_\lambda=E_\lambda\) on the positive cone, it follows that
\(
\widehat c_\lambda^m=c_\lambda^m.
\)

The Mountain Pass Theorem therefore provides a nontrivial critical point \(v_\lambda\) of \(\widehat E_\lambda\) such that
\(
\widehat E_\lambda(v_\lambda)=c_\lambda^m\ge\eta_\lambda>0.
\)
Testing \(\widehat E_\lambda'(v_\lambda)=0\) with \(v_\lambda^-=\max\{-v_\lambda,0\}\), we obtain
\(
0
=
\widehat E_\lambda'(v_\lambda)v_\lambda^-
=
-\|v_\lambda^-\|^2.
\)
Thus \(v_\lambda^-=0\), and hence \(v_\lambda\ge0\). Since \(c_\lambda^m>0\), \(v_\lambda\not\equiv0\). Therefore \(v_\lambda\) is a nonnegative nontrivial critical point of \(E_\lambda\), and Proposition~\ref{pos-sol} yields
\(
v_\lambda>0~\text{in }\Omega.
\)
Hence \(v_\lambda\in D\).
\end{proof}
}

{
\begin{remark}
    It is important to emphasize that the mountain pass geometry is established using the minimizer associated with $c_\lambda^+$. Up to this point, we have ensured the existence of such a minimizer for parameters in the interval $\lambda \in (\lambda^* - \zeta, \lambda_1(\Omega))$. {Accordingly, the mountain pass result is asserted only for $\lambda \in (\lambda^* - \zeta, \lambda_1(\Omega))$.}
\end{remark}}

\section{{ Asymptotic Behavior and Proofs of the Main Theorems}}\label{SB4}
This section is dedicated to the analysis of the asymptotic behavior of the solutions obtained through the aforementioned results. To this end, define 
\begin{eqnarray*}
    {\bar{\lambda}:=\inf\{\lambda>0:~\eqref{eq1}~\mbox{has a nontrivial weak solution}\}}.
\end{eqnarray*}
{Since $\mathcal N_\lambda=\emptyset$ for $\lambda<\lambda_*$, no nontrivial solution exists below $\lambda_*$, and hence $\lambda_*\leq\bar\lambda$. On the other hand, Proposition~\ref{sol-N+} provides a nontrivial solution for some $\lambda<\lambda^*$, so $\bar\lambda<\lambda^*$.} In the next results, we obtain a more precise {estimate}.

\begin{proposition}\label{right-continuation}
Assume that $2<q<p<2^*$ and that $(b_0)$ holds.
Let $\mu<\lambda_1(\Omega)$. Then the following assertions hold:
\begin{itemize}
\item[$(i)$] If Problem~\eqref{eq1} with parameter $\lambda=\mu$
admits a nontrivial weak solution, then it admits a nonnegative
nontrivial weak solution.

\item[$(ii)$] If Problem~\eqref{eq1} with parameter $\lambda=\mu$
admits a nonnegative nontrivial weak solution, then, for every
$\lambda\in(\mu,\lambda_1(\Omega))$, Problem~\eqref{eq1}
admits a nonnegative nontrivial weak solution.
\end{itemize}
\end{proposition}
\begin{proof}
Set $h(s)=s^{q-1}-s^{p-1}$ for $s\geq0$ and $M:=\max_{s\geq0}h(s)$. To prove $(i)$, let $w$ be a nontrivial weak solution for $\lambda=\mu$. {By the standard subcritical elliptic bootstrap, $w\in C^{1,\beta}(\overline\Omega)$ for some $\beta\in(0,1)$. By Kato's inequality, $|w|$ is a weak subsolution of}
$
-\Delta u=\mu u+b(x)h(u).
$
Since $\mu<\lambda_1(\Omega)$, the problem
$
-\Delta z-\mu z=M b(x)~\text{in }\Omega,~ z=0~\text{on }\partial\Omega,
$
has a positive solution $z$. For $K\geq1$,
\[
(-\Delta-\mu)(Kz)=KM b(x)\geq b(x)h(Kz),
\]
so $Kz$ is a supersolution. By the boundary point lemma and standard boundary estimates
(see, e.g., \cite{pucci}), there exists $K>0$ sufficiently
large such that
$
|w|\leq Kz\quad\text{in }\Omega.
$
Consequently, $Kz$ is a supersolution lying above the
subsolution $|w|$. The standard sub--supersolution theorem
(see, e.g., \cite{drabek}) therefore yields a nonnegative
nontrivial weak solution.

For $(ii)$, let $u_\mu\geq0$ be a nontrivial solution and let $\mu<\lambda<\lambda_1(\Omega)$. Then
\[
-\Delta u_\mu=\mu u_\mu+b(x)h(u_\mu)\leq \lambda u_\mu+b(x)h(u_\mu),
\]
so $u_\mu$ is a subsolution for the problem with parameter $\lambda$. Let $z_\lambda>0$ solve
\begin{equation}\label{quasi-super}
-\Delta z_\lambda-\lambda z_\lambda=M b(x)\quad\text{in }\Omega,\qquad z_\lambda=0\quad\text{on }\partial\Omega.
\end{equation}
As above, $Kz_\lambda$ is a supersolution for every $K\geq1$, and, {by the same boundary comparison, for $K$ sufficiently large one has $u_\mu\leq Kz_\lambda$. The sub--supersolution theorem \cite{drabek} therefore gives a nonnegative nontrivial solution for every $\mu<\lambda<\lambda_1(\Omega)$.}
\end{proof}
}

\begin{proposition}\label{lambda-bar}
Problem~\eqref{eq1} admits at least one positive solution when $\lambda=\bar\lambda$.
\end{proposition}
\begin{proof}

{
Let
$
\mathscr S:=\{\lambda>0:(P_\lambda)\text{admits a nontrivial weak solution}\}.
$
By definition, $\bar\lambda=\inf\mathscr S$.

If $\bar\lambda\in\mathscr S$, Proposition~\ref{right-continuation}$(i)$ gives a nonnegative nontrivial solution of $(P_{\bar\lambda})$. Proposition~\ref{pos-sol} then implies that this solution is positive.

Assume therefore, by contradiction, that
\(
\bar\lambda\notin\mathscr S.
\)
Since Proposition~\ref{sol-N+} gives solutions for parameters below $\lambda^*$, we have
\(
\bar\lambda<\lambda^*<\lambda_1(\Omega).
\)
Thus, there exists a sequence $(\lambda_k)\subset\mathscr S$ such that
\(
\bar\lambda<\lambda_k<\lambda^*,
~
\lambda_k\downarrow\bar\lambda.
\)
By Proposition~\ref{right-continuation}$(i)$, for each $k$ we may choose a nonnegative nontrivial solution $u_k$ of $(P_{\lambda_k})$.

Testing the equation with $u_k$, we obtain
\[
\|u_k\|^2-\lambda_k\|u_k\|_2^2
=
\int_\Omega b(x)(u_k^q-u_k^p)\,dx.
\]
Since $t^q-t^p$ is bounded above on $[0,\infty)$ and
\(
\lambda_k\leq\lambda^*<\lambda_1(\Omega),
\)
there exists $C>0$, independent of $k$, such that
\[
\left(
1-\frac{\lambda^*}{\lambda_1(\Omega)}
\right)\|u_k\|^2
\leq C.
\]
Hence $(u_k)$ is bounded in $H_0^1(\Omega)$.

On the other hand, the Nehari identity and the uniform bound
$\lambda_k\leq\lambda^*$ give, exactly as in Lemma~\ref{coercive}, a constant
$\delta>0$, independent of $k$, such that
\(
\|u_k\|\geq\delta.
\)
Therefore, up to a subsequence,
\(
u_k\rightharpoonup u
~\text{in }H_0^1(\Omega),
\) and
\(
u_k\to u
~\text{in }L^r(\Omega),
~
2\leq r<2^*.
\)

Testing the equations with $u_k-u$ and using the compact convergence of the subcritical terms, we obtain
\(
u_k\to u
~\text{strongly in }H_0^1(\Omega).
\)
Consequently,
\(
u\not\equiv0,
~
u\geq0.
\)
Passing to the limit in the weak formulation, we conclude that $u$ is a nontrivial weak solution of $(P_{\bar\lambda})$. Hence
\(
\bar\lambda\in\mathscr S,
\)
which contradicts the assumption $\bar\lambda\notin\mathscr S$.

Therefore $\bar\lambda\in\mathscr S$. Proposition~\ref{right-continuation}$(i)$ yields a nonnegative nontrivial solution $u_{\bar\lambda}$ of $(P_{\bar\lambda})$, and Proposition~\ref{pos-sol} gives
\(
u_{\bar\lambda}>0
~\text{in }\Omega.
\)
}
\end{proof}
{
\begin{proposition}\label{low-multiplicity}
Assume \(2<q<p<2^*\) and that \((b_0)\) holds. Then, for every
\(
\bar\lambda<\lambda\leq\lambda^*-\zeta,
\)
Problem~\eqref{eq1} admits at least two distinct positive solutions.
\end{proposition}

\begin{proof}
Let \(u_{\bar\lambda}>0\) be the solution of \((P_{\bar\lambda})\) given by Proposition~\ref{lambda-bar}. Since \(\lambda>\bar\lambda\), \(u_{\bar\lambda}\) is a strict subsolution of \((P_\lambda)\). On the other hand, as in Proposition~\ref{right-continuation}, for \(K>1\) sufficiently large, \( \overline u_\lambda:=Kz_\lambda \)
is a supersolution satisfying
\(
u_{\bar\lambda}\leq\overline u_\lambda.
\)

Set
\(
f_\lambda(x,t)
=
\lambda t+b(x)\bigl(t^{q-1}-t^{p-1}\bigr),
~t\geq0,
\)
and define the truncation
\[
\widetilde f_\lambda(x,t)
=
\begin{cases}
f_\lambda(x,u_{\bar\lambda}(x)),
& t<u_{\bar\lambda}(x),\\
f_\lambda(x,t),
& u_{\bar\lambda}(x)\leq t\leq\overline u_\lambda(x),\\
f_\lambda(x,\overline u_\lambda(x)),
& t>\overline u_\lambda(x).
\end{cases}
\]
Let \(\widetilde E_\lambda\) be the associated energy functional. Since
\(\widetilde f_\lambda\) is bounded, \(\widetilde E_\lambda\) is coercive and weakly lower semicontinuous, and therefore it has a global minimizer \(w_\lambda\).

The standard truncation argument gives
\(
u_{\bar\lambda}\leq w_\lambda\leq\overline u_\lambda
~\text{in }\Omega.
\)
Hence the truncation is inactive at \(w_\lambda\), and therefore \(w_\lambda\) is a weak solution of \((P_\lambda)\).

We next show that both inequalities are strict. Since \(u_{\bar\lambda}\) solves \((P_{\bar\lambda})\), setting
\(
\phi:=w_\lambda-u_{\bar\lambda}\geq0,
\)
we have
\[
-\Delta\phi
=
\lambda\phi
+
(\lambda-\bar\lambda)u_{\bar\lambda}
+
b(x)\bigl(h(w_\lambda)-h(u_{\bar\lambda})\bigr),
\]
where \(h(t)=t^{q-1}-t^{p-1}\). Since \(u_{\bar\lambda}\) and \(w_\lambda\) are bounded and \(h\in C^1([0,\infty))\), there exists \(C>0\) such that
\[
-\Delta\phi+C\phi
\geq
(\lambda-\bar\lambda)u_{\bar\lambda}>0
\qquad\text{in }\Omega.
\]
Thus, by the strong maximum principle and the boundary point lemma,
\(
u_{\bar\lambda}<w_\lambda
~\text{in }\Omega.
\)

Similarly, setting
\(
\psi:=\overline u_\lambda-w_\lambda\geq0,
\)
and recalling that \(\overline u_\lambda=Kz_\lambda\), we obtain
\[
-\Delta\psi-\lambda\psi
=
b(x)\bigl(KM-h(w_\lambda)\bigr)\geq0.
\]
Since \(K>1\), \(h(w_\lambda)\leq M\), and \(b>0\) a.e. on
\(\Omega\setminus\Omega_0\), the right-hand side is not identically zero. Hence the strong maximum principle and the boundary point lemma yield
\(
w_\lambda<\overline u_\lambda
~\text{in }\Omega.
\)

Therefore \(w_\lambda\) lies in the interior, with respect to the
\(C_0^1(\overline\Omega)\)-topology, of the order interval
\(
[u_{\bar\lambda},\overline u_\lambda].
\)
Consequently, there exists a \(C_0^1(\overline\Omega)\)-neighborhood of
\(w_\lambda\) in which the truncated functional \(\widetilde E_\lambda\)
coincides with \(E_\lambda\). Since \(w_\lambda\) is a global minimizer of
\(\widetilde E_\lambda\), it follows that \(w_\lambda\) is a local minimizer of
\(E_\lambda\) in \(C_0^1(\overline\Omega)\). By the \(H_0^1\)-versus-\(C^1\)
local minimizer result of Brezis and Nirenberg~\cite{BrezisNirenbergMin},
\(w_\lambda\) is also a local minimizer of \(E_\lambda\) in
\(H_0^1(\Omega)\).

Now consider
\[
\widehat E_\lambda(u)
=
\frac12\|u\|^2
-\frac{\lambda}{2}\|u^+\|_2^2
+\frac1p\|u^+\|_{p,b}^p
-\frac1q\|u^+\|_{q,b}^q.
\]
Since \(w_\lambda>0\) in \(\Omega\) and satisfies the boundary point property, it lies in the interior of the positive cone of \(C_0^1(\overline\Omega)\). Hence \(\widehat E_\lambda=E_\lambda\) in a \(C_0^1(\overline\Omega)\)-neighborhood of \(w_\lambda\), and therefore \(w_\lambda\) is a local minimizer of \(\widehat E_\lambda\) in \(C_0^1(\overline\Omega)\). By the \(H_0^1\)-versus-\(C^1\) local minimizer result of Brezis and Nirenberg~\cite{BrezisNirenbergMin}, \(w_\lambda\) is also a local minimizer of \(\widehat E_\lambda\) in \(H_0^1(\Omega)\).

On the other hand, \(0\) is a strict local minimizer of \(\widehat E_\lambda\), since \(\lambda<\lambda_1(\Omega)\), and \(\widehat E_\lambda\) satisfies the Palais--Smale condition. Therefore, by the Pucci--Serrin three critical point theorem~\cite{PucciSerrin}, \(\widehat E_\lambda\) possesses a critical point
\(
v_\lambda\notin\{0,w_\lambda\}.
\)
Testing \(\widehat E_\lambda'(v_\lambda)=0\) with \(v_\lambda^-\) gives
\(v_\lambda\geq0\). Therefore \(v_\lambda\) is a nonnegative nontrivial solution of \((P_\lambda)\), and Proposition~\ref{pos-sol} implies
\(
v_\lambda>0~\text{in }\Omega.
\)
Hence \(w_\lambda\) and \(v_\lambda\) are two distinct positive solutions.
\end{proof}
}

{Following \cite{KayeAb}, we will show that Problem~\eqref{eq1} has no nonnegative nontrivial solution at the level $\lambda_*$. Combined with the existence of a nonnegative solution at $\bar\lambda$, this yields $\lambda_*<\bar\lambda$ in the subsequent corollary.}

{
\begin{proposition}\label{lambda*-nonexist}
Suppose $2<q<p<2^*$ and that condition $(b_0)$ holds. Then Problem \eqref{eq1} has no nonnegative nontrivial solution for $\lambda=\lambda_*$.
\end{proposition}
\begin{proof}
Suppose, by contradiction, that $u\geq0$, $u\not\equiv0$, is a solution of Problem \eqref{eq1} for $\lambda=\lambda_*$. Then $u\in\mathcal{N}_{\lambda_*}$. If $u\notin D$, the nonlinear terms vanish and $u\in H_0^1(\Omega_0)$, so
\[
\lambda_*=\frac{\|u\|^2}{\|u\|_2^2}\geq\lambda_1(\Omega_0)>\lambda_*,
\]
a contradiction. Hence $u\in D$ and
\(
R_n(u)=\lambda_*.
\)
Since
$
\lambda_n(u)=\min_{t>0}R_n(tu)\leq R_n(u)=\lambda_*
$
and $\lambda_n(u)\geq\lambda_*$ by the definition of $\lambda_*$, we obtain
\(
\lambda_n(u)=R_n(u)=\lambda_*.
\)
Thus $t=1$ realizes the minimum of $t\mapsto R_n(tu)$. By the uniqueness of the minimum, $\tau(u)=1$, and hence $u\in\mathcal{N}_{\lambda_*}^0$. Moreover, $u$ is a minimizer of $\lambda_n$. Therefore, Proposition~\ref{ln} yields
\begin{equation}\label{II}
-2\Delta u=2\lambda_*u-b(x)(p|u|^{p-2}u-q|u|^{q-2}u).
\end{equation}
On the other hand, $u$ satisfies
\begin{equation}\label{I}
-\Delta u=\lambda_*u-b(x)(|u|^{p-2}u-|u|^{q-2}u).
\end{equation}
By Proposition~\ref{pos-sol}, $u>0$ in $\Omega$. Combining \eqref{I} and \eqref{II}, we obtain
$
b(x)\big((p-2)u^{p-1}-(q-2)u^{q-1}\big)=0.
$
{Since $b>0$ a.e. on $\Omega\setminus\Omega_0$ by $(b_0)$,
\[
u^{p-q}=\frac{q-2}{p-2},
\qquad\text{a.e. on }\Omega\setminus\Omega_0.
\]
The regularity used in Proposition~\ref{pos-sol} gives continuity up to $\partial\Omega$. Hence $u$ equals the same positive constant on $\Omega\setminus\overline{\Omega_0}$. Since $\Omega_0\subset\subset\Omega$, points of this set approach $\partial\Omega$; continuity would therefore force the same positive boundary value, contradicting $u=0$ on $\partial\Omega$.}
\end{proof}
}
\begin{cor}\label{comp-par}
$\lambda_*<\bar\lambda$.
\end{cor}
\begin{proof}
{
If $\bar\lambda=\lambda_*$, Proposition~\ref{lambda-bar} gives a positive solution of Problem~\eqref{eq1} at $\lambda=\lambda_*$, contradicting Proposition~\ref{lambda*-nonexist}. Hence
\(
\lambda_*<\bar\lambda.
\)
}
\end{proof}
Next, we will consider the asymptotic behavior of the solution on $\mathcal{N}_\lambda^+$.
{
\begin{proposition}\label{u-lamb-inf}
    Assume $2 < q < p < 2^*$ and that condition $(b_0)$ holds.
Let $u_\lambda \in \mathcal{N}_\lambda^+$ be the minimizing solution provided by Theorem \ref{t2} for $\lambda^*-\zeta<\lambda<\lambda_1(\Omega_0)$. Then, as $\lambda \to \lambda_1(\Omega_0)^-$, we have $\|u_\lambda\| \to \infty$.
\end{proposition}
\begin{proof}
    Assume, by contradiction, that there exists a sequence $(\hat\lambda_k)$ such that $\hat\lambda_k\nearrow \lambda_1(\Omega_0)$ and $(u_{\hat\lambda_k})$ is bounded in $H_0^1(\Omega)$.

From Theorem \ref{t2}, we have {$c_{\hat\lambda_k}^+ = E_{\hat\lambda_k}(u_k)$},
where $u_k:=u_{\hat\lambda_k} \in \mathcal{N}_{\hat\lambda_k}^+$. Furthermore, for each $k\geq 1$, $u_k$ is a solution of the problem 
    \begin{equation*}\label{eq1k}
		\left\{
            \begin{array}{rcll}
			-\Delta u  &=&  \hat\lambda_k u-b(x)(|u|^{p-2}u- |u|^{q-2}u),&\mbox{in}   \  \Omega,  \\[2ex]
			u&=&0,  &\mbox{on} \ \ \ \ \partial \Omega,
		\end{array}
		\right.
	\end{equation*}    
    for $\overline{\lambda}<\hat \lambda_k<\lambda_1(\Omega_0)$. Since  $(u_k)$ is a bounded sequence in $H_0^1(\Omega)$,   
    considering a subsequence if necessary, there exists $\tilde{u}\in H_0^1(\Omega)$ such that $u_k$ {converges weakly} to $\tilde{u}$ in $H_0^1(\Omega)$. Analogously to {Proposition} \ref{lambda-bar}, we conclude that $u_k\to \tilde{u}$ strongly in $H_0^1(\Omega)$. Consequently, $\tilde{u}$ is a solution of Problem \eqref{eq1} with $\lambda=\lambda_1(\Omega_0)$.

{Fix $\lambda_0\in(\lambda^*,\lambda_1(\Omega_0))$. By Theorem~\ref{t2}, $c_{\lambda_0}^+<0$. Since Proposition~\ref{cont-cplus} ensures that the map $\lambda\mapsto c_\lambda^+$ is strictly decreasing on $(\lambda^*-\zeta,\lambda_1(\Omega_0))$, for $k$ sufficiently large, $\hat\lambda_k\geq\lambda_0$ and therefore $
c_{\hat\lambda_k}^+\leq c_{\lambda_0}^+<0.
$
Hence
\[
E_{\lambda_1(\Omega_0)}(\tilde u)
=\lim_{k\to\infty}E_{\hat\lambda_k}(u_k)
=\lim_{k\to\infty}c_{\hat\lambda_k}^+
\leq c_{\lambda_0}^+<0.
\]
Thus, $\tilde u\not\equiv 0$. Moreover, since $u_k\geq 0$ and
$u_k\to\tilde u$ strongly in $H_0^1(\Omega)$, we have
$\tilde u\geq 0$ a.e. in $\Omega$. Proposition~\ref{pos-sol}
then yields $\tilde u>0$ in $\Omega$. On the other hand, from Theorem \ref{t4}, there is no positive solution to Problem \eqref{eq1} when $\lambda=\lambda_1(\Omega_0)$. Thus, $\|u_\lambda\|\to \infty$  when $\lambda\nearrow \lambda_1(\Omega_0)$.
}
\end{proof}
{
\begin{proposition}\label{profile-u-lambda}
Assume $2<q<p<2^*$ and that condition $(b_0)$ holds. Let $u_\lambda\in\mathcal N_\lambda^+$ be the minimizing solution provided by Theorem~\ref{t2} for
\(
\lambda^*-\zeta<\lambda<\lambda_1(\Omega_0).
\)
Then
\[
\frac{u_\lambda}{\|u_\lambda\|}
\longrightarrow
\frac{\varphi_1^0}{\|\varphi_1^0\|},
\quad\text{strongly in }H_0^1(\Omega),
\]
as $\lambda\to\lambda_1(\Omega_0)^-$.
\end{proposition}

\begin{proof}
Let $(\lambda_k)$ be any sequence such that
\(
\lambda_k\nearrow\lambda_1(\Omega_0)
\)
and set
\(
u_k:=u_{\lambda_k},
~
w_k:={u_k}/{\|u_k\|}.
\)
By Proposition~\ref{u-lamb-inf},
\(
\|u_k\|\to\infty.
\)
Since $\|w_k\|=1$, up to a subsequence, there exists
$w\in H_0^1(\Omega)$ such that
\(
w_k\rightharpoonup w
~\text{in }H_0^1(\Omega),
~
w_k\to w
~\text{in }L^r(\Omega),
~ 2\le r<2^*.
\)

Since $u_k\in\mathcal N_{\lambda_k}^+$, the Nehari identity yields
\[
\|u_k\|^2-\lambda_k\|u_k\|_2^2
=
\int_\Omega b(x)\bigl(|u_k|^q-|u_k|^p\bigr)\,dx.
\]
Since the function $t\mapsto t^q-t^p$ is bounded above on $[0,\infty)$, there exists $C>0$, independent of $k$, such that
$
\|u_k\|^2-\lambda_k\|u_k\|_2^2\le C.
$
Dividing by $\|u_k\|^2$ and using Proposition~\ref{u-lamb-inf}, we obtain
\(
1-\lambda_k\|w_k\|_2^2\to0.
\)
Hence
\(
\lambda_1(\Omega_0)\|w\|_2^2=1.
\)

On the other hand, dividing the Nehari identity by $\|u_k\|^p$, we obtain
\[
\|u_k\|^{2-p}
\bigl(1-\lambda_k\|w_k\|_2^2\bigr)
=
\|u_k\|^{q-p}\|w_k\|_{q,b}^q
-
\|w_k\|_{p,b}^p.
\]
Since $p>q>2$, $\|u_k\|\to\infty$, and $(w_k)$ is bounded in
$H_0^1(\Omega)$, it follows that
\(
\|w_k\|_{p,b}^p\to0.
\)
Since $w_k\to w$ strongly in $L^p(\Omega)$, we obtain
\(
\|w\|_{p,b}=0.
\)
Thus, by $(b_0)$,
\(
w=0
~\text{a.e. in }\Omega\setminus\Omega_0,
\)
and consequently $w\in H_0^1(\Omega_0)$.

Therefore, by the variational characterization of
$\lambda_1(\Omega_0)$,
\[
1
=
\lambda_1(\Omega_0)\|w\|_2^2
\le
\|w\|^2
\le
\liminf_{k\to\infty}\|w_k\|^2
=
1.
\]
Hence
$
\|w\|^2
=
\lambda_1(\Omega_0)\|w\|_2^2
=
1.
$
Thus $w$ is a first eigenfunction of $-\Delta$ in $\Omega_0$. Moreover, since $u_k\ge0$, we have $w\ge0$, and therefore
\(
w={\varphi_1^0}/{\|\varphi_1^0\|}.
\)
Finally, since $w_k\rightharpoonup w$ in $H_0^1(\Omega)$ and
\(
\|w_k\|=\|w\|=1,
\)
we conclude that
\(
w_k\to w
~\text{strongly in }H_0^1(\Omega).
\)

Since the limit is independent of the sequence $(\lambda_k)$, the whole family satisfies
\[
\frac{u_\lambda}{\|u_\lambda\|}
\longrightarrow
\frac{\varphi_1^0}{\|\varphi_1^0\|},
\quad\text{strongly in }H_0^1(\Omega),
\]
as $\lambda\to\lambda_1(\Omega_0)^-$.
\end{proof}
}
{
\begin{cor}\label{concentration-u-lambda}
Under the assumptions of Proposition~\ref{profile-u-lambda}, as
$\lambda\to\lambda_1(\Omega_0)^-$,
\[
\frac{\|u_\lambda\|_2^2}{\|u_\lambda\|^2}
\longrightarrow
\frac{1}{\lambda_1(\Omega_0)}
\]
and
\[
\frac{\|u_\lambda\|_{r,\Omega\setminus\Omega_0}}
{\|u_\lambda\|}
\longrightarrow0,
\qquad 2\le r<2^*.
\]
In particular,
\[
\frac{\|u_\lambda\|_{p,b}+\|u_\lambda\|_{q,b}}
{\|u_\lambda\|}
\longrightarrow0.
\]
\end{cor}
\begin{proof}
The result follows directly from Proposition~\ref{profile-u-lambda}, the compact embeddings
$H_0^1(\Omega)\hookrightarrow L^r(\Omega)$, $2\le r<2^*$, and the facts that
$\varphi_1^0=0$ a.e. in $\Omega\setminus\Omega_0$ and $b=0$ in $\Omega_0$.
\end{proof}
}
The following results play a crucial role in  establishing properties related to the asymptotic behavior of the mountain pass type solutions.
\begin{lemma}\label{limtlambda1}
    Assume that $2 < q < p < 2^*$ and condition $(b_0)$ holds. If $\lambda \in (\lambda_e(\varphi_1), \lambda_1(\Omega))$, then the function $\lambda \mapsto t_{\lambda,-}(\varphi_1)$ is differentiable and strictly decreasing. In particular, the sequence $(t_{\lambda,-}(\varphi_1))_\lambda$ is bounded.
\end{lemma}
\begin{proof}
Fix a parameter $\lambda' \in (\lambda_e(\varphi_1), \lambda_1(\Omega))$ and define the function $F: (\lambda_e(\varphi_1), \lambda_1(\Omega)) \times \mathbb{R} \to \mathbb{R}$ by
$$F(\lambda, t) = R_n(t\varphi_1) - \lambda.$$
It is straightforward to deduce that $F(\lambda', t_{\lambda',-}(\varphi_1)) = 0$ and, since $t_{\lambda',-}(\varphi_1)<\tau(\varphi_1)$,
$$
\frac{\partial F}{\partial t}(\lambda', t_{\lambda',-}(\varphi_1))
=
\left.\frac{d}{dt}R_n(t\varphi_1)\right|_{t=t_{\lambda',-}(\varphi_1)}
<0.
$$
By the Implicit Function Theorem, there exists $\varepsilon_0 > 0$ and a differentiable function $f: (\lambda' - \varepsilon_0, \lambda' + \varepsilon_0) \to \mathbb{R}$ such that $F(\lambda, f(\lambda)) = 0$,
$$
\left.\frac{d}{dt}R_n(t\varphi_1)\right|_{t=f(\lambda)} < 0
\quad\text{and}\quad
\frac{d}{d\lambda}f(\lambda)
=
\left(
\left.\frac{d}{dt}R_n(t\varphi_1)\right|_{t=f(\lambda)}
\right)^{-1}
<0,
$$
for all $\lambda \in (\lambda' - \varepsilon_0, \lambda' + \varepsilon_0)$.

Since $t_{\lambda,-}(\varphi_1)$ is the unique root of $R_n(t\varphi_1)=\lambda$ on the decreasing branch $0<t<\tau(\varphi_1)$, we have $f(\lambda)=t_{\lambda,-}(\varphi_1)$. As this conclusion holds for any $\lambda' \in (\lambda_e(\varphi_1), \lambda_1(\Omega))$, the lemma holds for all $\lambda \in (\lambda_e(\varphi_1), \lambda_1(\Omega))$. Moreover, $0 < t_{\lambda,-}(\varphi_1)$ and $t_{\lambda,-}(\varphi_1)$ is strictly decreasing with respect to $\lambda$. Consequently, $(t_{\lambda,-}(\varphi_1))_{\lambda}$ is bounded as $\lambda\to \lambda_1(\Omega)^-$.
\end{proof}

\begin{proposition}\label{ass-1}
Assume $2<q<p<2^*$ and that condition $(b_0)$ holds.
Let $v_\lambda$ be the mountain-pass solution
provided by Lemma~\ref{mp-solution}. Then
$
c_\lambda^m\longrightarrow0
~\text{as }\lambda\to\lambda_1(\Omega)^-.
$
\end{proposition}

\begin{proof}
For $\lambda>\lambda_e(\varphi_1)$, choose
$e_\lambda=t_{\lambda,+}(\varphi_1)\varphi_1$.
Then the path
$
\gamma_\lambda(r)
=r\,t_{\lambda,+}(\varphi_1)\varphi_1,
~ r\in[0,1],
$
belongs to $\Gamma_\lambda$. Consequently,
\[
\begin{aligned}
0<c_\lambda^m
&\le \max_{r\in[0,1]}E_\lambda(\gamma_\lambda(r))=E_\lambda(t_{\lambda,-}(\varphi_1)\varphi_1)\\
&\le \frac{q-2}{2q}
t_{\lambda,-}^2(\varphi_1)
\bigl(\|\varphi_1\|^2
-\lambda\|\varphi_1\|_2^2\bigr)
=\frac{q-2}{2q}
t_{\lambda,-}^2(\varphi_1)
\bigl(\lambda_1(\Omega)-\lambda\bigr)
\|\varphi_1\|_2^2.
\end{aligned}
\]
By Lemma~\ref{limtlambda1},
$t_{\lambda,-}(\varphi_1)$ remains bounded
as $\lambda\to\lambda_1(\Omega)^-$.
Therefore, $c_\lambda^m\to0$.
\end{proof}

\noindent
We are now in a position to complete the proofs of our main theorems
using the results established above.

\noindent\textbf{Proof of Theorem \ref{t2}:} {By Proposition~\ref{lambda-bar} and Proposition~\ref{right-continuation}$(ii)$, Problem~\eqref{eq1} admits a non-negative nontrivial solution for every $\lambda\in(\bar\lambda,\lambda_1(\Omega))$. 

Since $\bar\lambda<\lambda^*$, we may decrease
$\zeta>0$, if necessary, so that
$\bar\lambda<\lambda^*-\zeta$.
By Proposition~\ref{sol-N+},
Lemma~\ref{min-sol}, and
Proposition~\ref{belo-D}, there exists a
minimizing solution
$u_\lambda\in\mathcal N_\lambda^+\cap D$
for every
$\lambda\in(\lambda^*-\zeta,\lambda_1(\Omega_0))$.
Moreover, since $E_\lambda(|u|)=E_\lambda(u)$
and both the Nehari identity and the sign of
the second derivative of the fibering map
are preserved under the replacement of $u$
by $|u|$, this solution may be chosen
non-negative. Proposition~\ref{pos-sol}
then implies that $u_\lambda$ is positive
in $\Omega$.

In particular, the variational construction provides positive solutions for every $\lambda\in[\lambda_1(\Omega),\lambda_1(\Omega_0))$. Combined with the continuation argument, this establishes the existence of a positive solution for every $\lambda\in(\bar\lambda,\lambda_1(\Omega_0))$.}

We now establish the energy assertions.
For $\lambda^*-\zeta<\lambda<\lambda^*$,
let $u_\lambda\in\mathcal N_\lambda^+\cap D$
be the minimizing solution.
Since $u_\lambda\in D$, the definition of
$\lambda^*$ implies that
$
R_e(u_\lambda)\geq\lambda^*.
$
Consequently, using
$
E_\lambda(u_\lambda)
=1/2\|u_\lambda\|_2^2
\bigl(R_e(u_\lambda)-\lambda\bigr),
$
we obtain $E_\lambda(u_\lambda)>0$.

At $\lambda=\lambda^*$,
Proposition~\ref{le} provides a solution
$u_{\lambda^*}$ such that
$E_{\lambda^*}(u_{\lambda^*})=0$.

Finally, for
$\lambda^*<\lambda<\lambda_1(\Omega_0)$,
the identity
$\lambda_e(u_{\lambda^*})=\lambda^*$
and Corollary~\ref{energ-signal},
items $(3)$--$(4)$, yield
$
E_\lambda\bigl(
t_{\lambda,+}(u_{\lambda^*})u_{\lambda^*}
\bigr)<0.
$
By the variational characterization of
$c_\lambda^+$, we have
$
E_\lambda(u_\lambda)=c_\lambda^+
\leq
E_\lambda\bigl(
t_{\lambda,+}(u_{\lambda^*})u_{\lambda^*}
\bigr)<0.
$
This establishes items $(i)$--$(iii)$.

Item $(iv)$ follows from
Proposition~\ref{u-lamb-inf},
item $(v)$ from Proposition~\ref{cont-cplus},
and item $(vi)$ from
Proposition~\ref{profile-u-lambda}.
\hfill$\square$

\

\noindent\textbf{Proof of Theorem~\ref{t1}.}
For $\bar\lambda<\lambda\leq\lambda^*-\zeta$,
the result follows from
Proposition~\ref{low-multiplicity}.

Now let
$\lambda^*-\zeta<\lambda<\lambda_1(\Omega)$.
By Theorem~\ref{t2} and
Lemma~\ref{mp-solution},
Problem~\eqref{eq1} admits positive solutions
$u_\lambda\in\mathcal N_\lambda^+$
and $v_\lambda$, respectively.

If $\lambda<\lambda^*$,
Lemma~\ref{geometry} yields
$
E_\lambda(u_\lambda)=c_\lambda^+
<
\eta_\lambda
\leq c_\lambda^m
=E_\lambda(v_\lambda).
$
On the other hand, if $\lambda\geq\lambda^*$,
$
E_\lambda(u_\lambda)=c_\lambda^+
\leq 0<c_\lambda^m
=E_\lambda(v_\lambda).
$
Thus, $u_\lambda\neq v_\lambda$,
and Problem~\eqref{eq1} admits at least
two distinct positive solutions.

Finally, the asymptotic behavior of
$c_\lambda^m$ follows from
Proposition~\ref{ass-1}.
\hfill$\square$

\

\noindent\textbf{Proof of Theorem \ref{t3}:} The result follows directly from Propositions \ref{lambda-bar} and \ref{pos-sol}.
\hfill$\square$

\ 

\noindent\textbf{Proof of Theorem~\ref{t4}.}
If $\lambda<\lambda_*$, then
$\mathcal N_\lambda=\emptyset$,
and hence Problem~\eqref{eq1}
admits no nontrivial solution.
Moreover, for
$\lambda_*\leq\lambda<\bar\lambda$,
the nonexistence of nontrivial solutions
follows directly from the definition
of $\bar\lambda$.

{Now, following an argument similar to that of Cardoso et al.~\cite{furtadomaia}, suppose, by contradiction, that $\lambda\geq\lambda_1(\Omega_0)$ and that $u>0$ is a solution of \eqref{eq1}. Extending $\varphi_1^0$ by zero outside $\Omega_0$ and using it as a test function, we obtain
\begin{equation}\label{u sol}
\int_{\Omega_0}\nabla u\cdot\nabla\varphi_1^0\,dx
=
\lambda\int_{\Omega_0}u\varphi_1^0\,dx,
\end{equation}
since $b=0$ a.e. in $\Omega_0$. 

On the other hand, integration by parts
over $\Omega_0$ yields
\[
\int_{\Omega_0}\nabla u\cdot\nabla\varphi_1^0\,dx
=
\lambda_1(\Omega_0)\int_{\Omega_0}u\varphi_1^0\,dx
+
\int_{\partial\Omega_0}u\,\partial_\eta\varphi_1^0\,d\sigma,
\]
where $\eta$ denotes the outward unit normal
to $\partial\Omega_0$.

Since $\overline{\Omega}_0\subset \subset\Omega$, the continuity and positivity of $u$ imply that $u>0$ on $\overline{\Omega}_0$, while Hopf's lemma gives
\(
\partial_\eta\varphi_1^0<0
~\text{on }\partial\Omega_0.
\)
Consequently,
\[
\int_{\Omega_0}\nabla u\cdot\nabla\varphi_1^0\,dx
<
\lambda_1(\Omega_0)\int_{\Omega_0}u\varphi_1^0\,dx.
\]
Combining this inequality with
\eqref{u sol}, we obtain
\[
(\lambda-\lambda_1(\Omega_0))
\int_{\Omega_0}u\varphi_1^0\,dx<0,
\]
which contradicts $\lambda\geq\lambda_1(\Omega_0)$.
\hfill$\square$}

\appendix
\section{}
\renewcommand{\theequation}{\thesection\arabic{equation}}  
\setcounter{equation}{0} 
\begin{definition}
    Let  $1 \leq r < \infty$ and \(b: \Omega \to [0, \infty)\) a measurable function called the weight function.  We define the space $$L_{r,b}(\Omega)
     = \left\{ u : \Omega \to \mathbb{R}; \|u\|_{L^r_b} < \infty \right\},
    \text{ where }
    ||u||_{L_{r,b}(\Omega)} :=||u||_{r,b}= \left( \int_\Omega b(x)|u|^r dx\right)^{1/r}.$$ 
\end{definition}
Under assumption $(b_0)$, for every $1\le r<\infty$,
we have
\[
D=\{u\in H_0^1(\Omega)\setminus\{0\}:\|u\|_{r,b}>0\}.
\]
Indeed, this follows from the fact that $b=0$ a.e. on $\Omega_0$ and $b>0$ a.e. on $\Omega\setminus\Omega_0$.}

\begin{remark}\label{des-r-b}
Let $u\in H_0^1(\Omega)$.
\begin{enumerate}
\item[$(1)$] If $1\le r<2^*$, then
\begin{equation}\label{imersao1}
\|u\|_{r,b}\leq\frac{\|b\|_\infty^{1/r}}{S_r}\|u\|.
\end{equation}
\item[$(2)$] If $1\le r<s<2^*$, then
\begin{equation}\label{desigualdade}
\|u\|_{r,b}^r\leq |\Omega|^{\frac{s-r}{s}}\|b\|_\infty^{\frac{s-r}{s}}\|u\|_{s,b}^r.
\end{equation}
\end{enumerate}
\end{remark}

\begin{lemma}
Assume that condition $(b_0)$ holds.
Then $D$ is open in $H_0^1(\Omega)\setminus\{0\}$
and $tD=D$ for every $t>0$.
\end{lemma}

\begin{proof}
By \eqref{imersao1}, the map
$u\mapsto\|u\|_{r,b}$ is continuous on
$H_0^1(\Omega)$. Consequently,
\[
D=\{u\in H_0^1(\Omega)\setminus\{0\}:
\|u\|_{r,b}>0\}
\]
is open in $H_0^1(\Omega)\setminus\{0\}$.

Moreover, since
$\|tu\|_{r,b}=t\|u\|_{r,b}$ for every $t>0$,
we have $u\in D$ if and only if $tu\in D$.
Hence $tD=D$.
\end{proof}


	\subsection*{Acknowledgements}
		\subsection*{Author Contribution}	{All authors contributed equally to this work.} {All authors also reviewed the manuscript in its present version.}
	
	\subsection*{Funding} 
	
	\subsection*{Availability of data and material} For the present work the availability of data and material is not applicable.
	
	\subsection*{Ethics approval and consent to participate} This item is not applicable for the present work.
	
\subsection*{Consent for publication} All authors read and approved the final version for the present manuscript.
	
	\subsection*{Declarations}
	
	\subsection*{Conflict of Interest}
	The authors declare that they have no conflict of interest.

\end{document}